\documentclass{jT}

\usepackage{amssymb,amsmath,latexsym,amsthm}
\usepackage[english]{babel}

\usepackage{euler, epic,eepic,latexsym, amssymb, amscd, amsfonts, xypic, floatflt}
\usepackage[OT2,T1]{fontenc}
\DeclareSymbolFont{cyrletters}{OT2}{wncyr}{m}{n}
\DeclareMathSymbol{\Sha}{\mathalpha}{cyrletters}{"58}

\input xy
\xyoption{all}

\newtheorem{theorem}{Theorem}[section]
\newtheorem{lemma}{Lemma}[section]
\newtheorem{proposition}{Proposition}[section]
\newtheorem{corollary}{Corollary}[section]

\theoremstyle{definition}
\newtheorem{definition}{Definition}

\numberwithin{equation}{section}

\newcommand{\Z}{\mathbb{Z}}

\newcommand{\Q}{\mathbb{Q}}
\newcommand{\Qbar}{\overline{\mathbb{Q}}}
\newcommand{\kbar}{\overline{k}}
\newcommand{\R}{\mathbb{R}}
\newcommand{\C}{\mathbb{C}}
\newcommand{\F}{\mathbb{F}}
\newcommand{\G}{\mathbb{G}}

\newcommand{\proj}{\mathbb P}

\newcommand{\tA}{\tilde{A}}
\newcommand{\tX}{\tilde{X}}
\newcommand{\tY}{\tilde{Y}}

\newcommand{\cF}{{\mathcal{F}}}
\newcommand{\cG}{{\mathcal{G}}}
\newcommand{\oh}{{\mathcal{O}}}

\renewcommand{\cL}{{\mathcal{L}}}
\newcommand{\cm}{{\mathcal{M}}}
\newcommand{\ct}{{\mathcal{T}}}

\newcommand{\cP}{{\mathcal{P}}}

\newcommand{\fT}{\mathfrak{T}}

\newcommand{\De}{\Delta}
\newcommand{\si}{\sigma}

\newcommand{\chr}{\operatorname{char}}

\newcommand{\Gal}{\operatorname{Gal}}
\newcommand{\Ad}{\operatorname{Ad}}
\newcommand{\uGal}{\underline{\Gal}}
\newcommand{\Hom}{\operatorname{Hom}}

\newcommand{\Maps}{\operatorname{Maps}}

\newcommand{\Pic}{\operatorname{Pic}}
\newcommand{\Spec}{\operatorname{Spec}}

\newcommand{\Aut}{\operatorname{Aut}}

\newcommand{\upi}{\underline{{\pi}_1}}
\newcommand{\red}{\rm{red}}

\begin{document}

\title[Universal covering spaces and fundamental groups in algebraic geometry]{Universal covering spaces and fundamental groups in algebraic geometry as schemes}


\author{\sc Ravi Vakil}
\address{Ravi Vakil \\
Department of Mathematics, Stanford University \\
Stanford CA USA 94305}
\email{vakil@math.stanford.edu}

\author{\sc Kirsten Wickelgren}
\address{Kirsten Wickelgren \\
Dept. of Mathematics, Harvard University \\
One Oxford St. \\
Cambridge MA USA 02138}
\email{kwickelg@math.harvard.edu}

\maketitle

\begin{resume}
En topologie, les notions du groupe fondamental et du rev\^etement universel sont lie\'es l'une \`a l'autre. En suivant l'exemple topologique, on construit une sch\'ema en groupes fondamentals  d'un rev\^etement universel, qui sont tous les deux des sch\'emas. Une fibre g\'eom\'etrique de la sch\'ema en groupes fondamentals est hom\'eomorphe au groupe fondamental \'etale. Ces constructions s'appliquent \`a toute sch\'ema quasicompact et quasis\'epar\'e. Avec des m\'ethodes et des hypoth\`eses differentes, cette sch\'ema en groupes fondamentals a d\'ej\`a \'et\'e construite par Deligne. 
\end{resume}

\begin{abstr}
In topology, the notions of the fundamental group and the universal cover are closely intertwined. By importing usual notions from topology into the algebraic and arithmetic setting, we construct a fundamental group family from a universal cover, both of which are schemes. A geometric fiber of the fundamental group
  family (as a topological group) is canonically the \'etale
  fundamental group. The constructions apply to all connected quasicompact quasiseparated
  schemes. With different methods and hypotheses, this fundamental group family was already constructed by Deligne.  \end{abstr}

\bigskip
\section{Introduction}

This paper takes certain natural topological constructions into the
algebraic and arithmetic setting. Primarily, we refer to the following: for a
sufficiently nice topological space $X$, the fundamental group
$\pi_1^{top}(X,x)$ varies continuously as $x$ varies. Thus, there is a
family of pointed fundamental groups, which we denote
$\xymatrix{\upi^{top}(X) \ar[r] &  X}$, whose fibers are canonically
$\pi_1^{top}(X,x)$.  The space  $\upi^{top}(X)$ is a group object among covering spaces.  We call it the {\em fundamental group family}. (It is also the isotropy group of the fundamental groupoid, as well as the adjoint bundle of the universal cover $\tX \rightarrow X$ viewed as a principal $\Aut(\tX/X)$-bundle, but both of these are awkwardly long to be used as names.) This paper repeats this process in the setting of algebraic geometry: for any connected quasicompact quasiseparated scheme $X$, we construct a group scheme $\upi(X) \rightarrow X$ whose fibers over geometric points are Grothendieck's \'etale fundamental group $\pi_1(X,x)$. This group scheme has already been constructed by Deligne  in \cite{deligne} for Noetherian schemes. The method and goals of the present paper are different from Deligne's, however, and we hope of interest. (See \S
\ref{relationtoearlierwork} for a further discussion of the fundamental group family in \cite{deligne}.)

The motivation for gluing together the $\pi_1(X,x)$ (which are
individually topological groups) into a group scheme requires some
explanation. We wish to study the question: what is a ``loop up to
homotopy'' on a scheme? Grothendieck's construction of the \'etale
fundamental group gives the beautiful perspective that loops up to
homotopy are what classify covering spaces. Although a map from the
circle to a scheme and the equivalence class of such a map up to
homotopy are problematic to define, \cite{sga1} defines the
fundamental group by first defining a covering space to be a finite
\'etale map, and then defining $\pi_1$ as the group classifying such
covering spaces. As finite \'etale maps of complex varieties are
equivalent to {\em finite} topological covering spaces, this
definition raises the question: why have we restricted to finite
covering spaces? There are at least two answers to this question,
neither of which is new: the first is that the
covering spaces of infinite degree may not be algebraic; it is the
finite topological covering spaces of a complex analytic space
corresponding to a variety that themselves correspond to
varieties. The second is that Grothendieck's \'etale $\pi_1$ 
classifies more than finite covers.  It classifies inverse limits of
finite \'etale covering spaces  \cite[Exp.~V.5, e.g., Prop.\ 5.2]{sga1}. These inverse limits are the profinite-\'etale
covering spaces we discuss in this paper (see Definition
\ref{d:pfe}). Grothendieck's enlarged fundamental group
\cite[Exp.~X.6]{sga3} even classifies some infinite covering spaces
that are not profinite-\'etale.

In topology, a covering space is defined to be a map which is locally trivial in the sense that it is locally of the form $\coprod U \rightarrow U$. We have the heuristic picture that to form a locally trivial space, you take a trivial space $\coprod U \rightarrow U$ and every time you go around a loop, you decide how to glue the trivial space to itself. (This heuristic picture is formalized by the theory of descent.) This leads to the notion that what the group of loops up to homotopy {\em should} classify are the locally trivial spaces. It becomes natural to ask: to what extent are finite \'etale or profinite-\'etale covering spaces locally trivial?\footnote{The same question should be asked for the covering spaces implicit in Grothendieck's enlarged fundamental group; we do not do this in this paper.} This is a substitute for the question: to what extent is \'etale $\pi_1$ the group of ``loops up to homotopy'' of a scheme?

The answer for finite \'etale maps is
well-known. (Finite \'etale maps are finite \'etale locally $\coprod_S
U \rightarrow U$ for $S$ a finite set. See \cite[Prop. 5.2.9]{sz} for a particularly enlightening exposition.) For profinite-\'etale maps, we
introduce the notion of Yoneda triviality and compare it to the notion
that a trivial map is a map of the form $\coprod U \rightarrow U$ (see
Definition \ref{Yoneda_trivial_def} and Proposition
\ref{Yoneda_triv_vs_disjoint_union}). Although a profinite-\'etale
morphism is locally Yoneda trivial (Corollary
\ref{pfe_cov_spaces_pfe_loc_Yon_triv}), locally Yoneda trivial
morphisms need not be profinite-\'etale. Indeed, the property of being
profinite-\'etale is not Zariski-local on the base (see Warning
\ref{warning}(b)). Since the \'etale fundamental group, which
classifies profinite-\'etale spaces, is obviously useful, but there
are other locally trivial spaces, this suggests that there are
different sorts of fundamental groups, each approximating ``loops up
to homotopy,'' by classifying some notion of a covering space, where a
covering space is some restricted class of locally trivial
spaces.\footnote{Note that the notion of a ``locally trivial space''
  is composed of the notion of ``locally'' and the notion of a
  ``trivial space.'' The idea of changing the notion of ``locally'' is
  thoroughly developed in the theory of Grothendieck topologies. Here,
  we are also interested in different notions of ``trivial.''} (Also
see \S \ref{K(pi,1)_elliptic_curve}.) See \cite{Morel} for the construction of the $\mathbb{A}^1$-fundamental group.

Returning to the motivation for constructing the fundamental group
family, it is not guaranteed that the object which classifies some
particular notion of covering space is a group; the \'etale
fundamental group is a topological group; and work of Nori \cite{nori}
shows that scheme structure can be necessary. (Nori's fundamental
group scheme is discussed in more detail in \S
\ref{relationtoearlierwork}.) However, a fiber of the fundamental
group family of \S \ref{s:FG} should classify covering spaces, and
indeed does in the case we deal with in this paper, where ``covering
space'' means profinite-\'etale morphism (see Theorem
\ref{etalefglink}).

More concretely, consider the following procedure: (1) define {\em
  trivial covering space}. (2) Define {\em covering space}. (3) Find a
large class of schemes which admit a simply connected covering space,
where a {\em simply connected} scheme is a scheme whose covering
spaces are all trivial. (4) Use (3) and the adjoint bundle construction
described in \S \ref{s:FG} to produce a fundamental group family. This
fundamental group family should be a group scheme over the base
classifying the covering spaces of (2).

We carry out this procedure with ``trivial covering space'' defined to
mean a Yoneda trivial profinite-\'etale morphism, and ``covering
space'' defined to mean a profinite-\'etale morphism. Then, for any
connected, quasicompact and quasiseparated scheme, there is a universal covering space
(see Proposition \ref{UCexists}), and the topological group underlying
the geometric fibers of the corresponding fundamental group family are the
pointed \'etale fundamental groups (see Theorem \ref{etalefglink}). In
particular, the topology on the \'etale fundamental group is the
Zariski topology on the geometric fibers of the fundamental group family.  Motivation for this
is the exercise that $\Spec \Qbar \otimes_{\Q} \Qbar$ (with the
Zariski topology) is homeomorphic to $\Gal (\Qbar / \Q)$ (with the
profinite
topology). 
We work through these ideas in a number of explicit examples.

\subsection{Relation to earlier work} \label{relationtoearlierwork} 
The fundamental group family of \S \ref{s:FG} can be recovered from Deligne's profinite fundamental group \cite[\S 10]{deligne}. For $X$ Noetherian, Deligne constructs a profinite lisse sheaf $P$ on $X \times X$, called the profinite fundamental group \cite[\S 10.17]{deligne}. $P$ pulled back under the diagonal map $X \rightarrow X \times X$, as in \cite[\S 10.5]{deligne}, gives the fundamental group family, via the equivalence between profinite lisse sheaves and pro-objects of the category of finite \'etale $X$ schemes. Deligne's constructions and goals are different from ours. For instance, he is not concerned with the fact that a sheaf can be locally profinite lisse without being so globally. His assumption that $X$ is Noetherian, however, is not necessary for \cite[\S 10.17]{deligne}; his construction should work for all quasicompact quasiseparated $X$, and in particular, our construction does not work in greater generality.

Over a
field $k$ and subject to additional hypotheses, other fundamental group
schemes have previously been constructed as well. Work of
Nori \cite{nori_compositio, nori} develops a pointed fundamental group scheme
which classifies principal $G$-bundles for $G$ a finite group scheme
over $k$, under additional hypotheses, including that the base scheme is connected,
reduced, and equipped with a rational point. The scheme structure is
necessary for this classification. Furthermore, Nori's fundamental group scheme has an associated
universal cover \cite[p.~84-85,~Def1~Prop2]{nori}. We expect that Nori's universal
cover admits a fundamental group family as in \S \ref{s:FG} whose
fiber over the given $k$-rational point is Nori's fundamental group
scheme. In particular, Nori's universal cover should not be the
universal cover of Proposition \ref{UCexists}.  We suspect that it is
the inverse limit of pointed principal $G$-bundles, where $G$ is a
finite group scheme over $k$, and that after base change to $\kbar$ and passage to the maximal pro-\'etale quotient, Nori's universal cover becomes the universal cover of Proposition \ref{UCexists}. We
have not verified these statements.

 Esnault and Hai \cite{eh}
define a variant of Nori's fundamental group scheme for a smooth
scheme $X$ over a characteristic $0$ field $k$, where $k$ is the field
of constants of $X$. The goals of Nori and Esnault-Hai are quite different from those of this paper.  For example, Esnault and Hai define a linear structure 
on Grothendieck's \'etale fundamental group of a scheme defined over a characteristic $0$ field.  They apply their results to study Grothendieck's
section conjecture.

The idea of changing the notion of ``covering space'' to recover the
classification of covering spaces by a fundamental group has appeared
earlier in topology. For example, Biss uses a fundamental group
equipped with a topology to classify ``rigid covering bundles'' over
some non semi-locally simply connected spaces (such as the Hawaiian
earring) \cite{biss1, biss2}, where the usual topological theory of
covering spaces is not valid. Moreover, ``rigid covering bundles,''
which are defined as Serre fibrations whose fibers have no
non-constant paths, are analogous to fiber bundles with totally
disconnected fiber. In the context of this paper, such a fiber bundle
should be viewed as a locally trivial space, where ``trivial'' is
defined to mean $U \times F \rightarrow U$, where $F$ is a totally
disconnected topological space. In earlier work,
Morgan and Morrison also considered the question of the ``right'' fundamental group of the Hawaiian earring and similar spaces in response
to a question of Eilenberg's, and with similar conclusions \cite{MM}.

Although we circumvent Noetherian hypotheses here, this is not new. Lenstra develops the theory of the \'etale fundamental group for any connected base scheme in his course on the Galois theory of schemes \cite{lenstra}. The \'etale fundamental group is only shown to classify finite (as opposed to profinite) \'etale covering spaces in this generality. He also stressed the importance of the local triviality of finite \'etale covering spaces, and drew close analogies to the topological theory.  See \cite{sz} for another nice exposition. 

The existence of the universal cover of Proposition \ref{UCexists} is also well-known to experts, but we include a proof in the required generality for completeness. 

The universal cover of a variety (in the sense of \S \ref{s:UC}) is
not in general a variety.  It is the algebraic analogue of a {\em
  solenoid} (see for example Dennis Sullivan's \cite{sullivan}), and
perhaps profinite-\'etale covering spaces of varieties deserve this
name as well.   The notion of solenoid in an algebraic context arose earlier, 
see \cite[p.~89]{artintate}.  Solenoids are examples of finite-dimensional
proalgebraic varieties in the sense of Piatetski-Shapiro and
Shafarevich, see \cite[\S 4]{Shaf}.  (Caution: Prop.\ 2 of \cite[\S
4]{Shaf} appears to be contradicted by Warning \ref{warning}(b).)
The notion essentially appears much earlier in Serre's \cite{serreGP}.

\subsection{Conventions}  As usual, fpqc means faithfully flat and
quasicompact, qcqs means quasicompact and quasiseparated, and $K^s$ is
the separable closure of $K$.  The phrase ``profinite-\'etale''
appears in the literature, but it is not clear to us that there is a
consistent definition, so to prevent confusion, we define it in
Definition~\ref{d:pfe}.  Warning: other definitions (such as the one
implicit in \cite{Shaf}) are different from ours, and would lead to a
{\em different} universal cover and fundamental group scheme.

\subsection{Acknowledgments}  We thank A. Blumberg, G. Carlsson,
R. Cohen, B. Conrad, T. Ekedahl, H. Esnault, J.
Ellenberg, K. Gruher, F. Grunewald, J. Hall, D. Harbater, S. Kerckhoff, M. Lieblich,
R. Lipshitz, C. McMullen, M. Nori, A. Ogus, M. Olsson, B.
Osserman, J. Rabinoff, D. Rydh,  S. Schr\"oer, T. Szamuely and C. Xu for many
helpful comments.  We thank R. Treger for pointing out \cite{Shaf}
to us.

The first author was
partially supported by NSF grant DMS-0801196. The second author was supported by an NSF Graduate Research Fellowship and a Stanford Graduate Fellowship.

\section{From topology to algebraic geometry, via a
``right'' notion of covering
space}

\label{s:coveringspace}
We now try to import the topological perspective into algebraic geometry in order to choose as parsimonious as possible a list of definitions, from which constructions and arguments are natural and straightforward.  Our endpoint is essentially the usual one; but we hope the reader may find the derivation direct.

\begin{definition} \label{Yoneda_trivial_def}A map of schemes $f: Y \rightarrow X$ is {\em Yoneda trivial} if $f$ admits a set of sections $S$ such that for each connected scheme $Z$, the natural map $$\Maps(Z,X) \times S \rightarrow \Maps(Z,Y) $$ is a bijection. Here, $\Maps(-,-)$ denotes the set of scheme morphisms. $S$ is called the {\em distinguished set of sections of $f$}.\end{definition}

The name ``Yoneda trivial'' comes from Yoneda's lemma, which controls $Y$ by the morphisms to $Y$; $Y$ is trivial over $X$ in the sense that maps to $Y$ from connected schemes are controlled by maps to $X$.

Note that if $X$ is connected, the distinguished sections must be the entire set of sections.

A trivial topological covering space is a map of topological spaces of the form $\coprod U \rightarrow U.$ We compare Yoneda trivial morphisms to morphisms of the form $\coprod X \rightarrow X$.

\begin{proposition}
{Let $X$ be a connected scheme. Then $\coprod X \rightarrow X$ is Yoneda trivial, 
where the coproduct is over any set. If $f: Y \rightarrow X$ is Yoneda trivial and the underlying topological space of $Y$ is a disjoint union of connected components (or if $Y$ is locally Noetherian), then $f$ is of the form $\coprod_S X \rightarrow X$ for some set $S$.}\label{Yoneda_triv_vs_disjoint_union}\end{proposition}

\begin{proof}
The first statement is obvious. For the second statement, we have $Y = \coprod_{c \in C} Y_c$ with $Y_c$ connected. Since $f$ is Yoneda trivial, the inclusion $Y_c \hookrightarrow Y$ factors through a distinguished section. It follows that $f: Y_c \rightarrow X$
is an isomorphism.
\end{proof}

The distinguished sections $S$ of a Yoneda trivial morphism $f: Y \rightarrow X$ can be given the structure of a topological space:  let $\fT$ denote the forgetful functor from schemes to topological spaces. It follows easily from the definition that Yoneda trivial morphisms induce isomorphisms on the residue fields of points, and therefore that the distinguished set of sections is in bijection with any fiber of $\fT(f): \fT(Y) \rightarrow \fT(X).$ In particular, $S$ is a subset of $\Maps_{cts}(\fT(X),\fT(Y))$, the continuous maps from $\fT(X)$ to $\fT(Y).$ Give $\Maps_{cts}(\fT(X),\fT(Y))$ the topology of pointwise convergence and give $S$ the subspace topology. 

\begin{definition} A morphism of schemes $f: Y \rightarrow X$ is {\em
  profinite-\'etale} if $Y = \underline{\Spec} \mathcal{A}$, where $\mathcal{A}$
is a colimit of quasi-coherent algebras, each corresponding to a finite \'etale morphism. Thus $f$ is
an inverse limit of finite \'etale morphisms.\label{d:pfe}\end{definition}

\begin{definition} A {\em covering space} is a profinite-\'etale
  morphism.\label{d:covering_space}
\end{definition}

We sometimes say (redundantly) {\em profinite-\'etale covering space}. (This redundancy comes from the point of view that there are other interesting notions of covering space.) 

Profinite-\'etale covering spaces are clearly stable under pull-back. Note that a profinite-\'etale covering space of a qcqs scheme is qcqs. By  \cite[\S 8 Th.\ 8.8.2]{ega4iii},  \cite[\S 8 Th.\ 8.10.5]{ega4iii}, and  \cite[\S 17 Cor.\ 17.7.3]{ega4iv} profinite-\'etale covering spaces of qcqs are closed under composition.

\subsection{Warnings}  \label{warning}
 (a) Although a
profinite-\'etale morphism is integral, flat, and formally unramified, the
converse need not hold.  For example, let $p$ be a prime, $X = \Spec
\F_p (t)$, and $$Y =\Spec \F_p(t^{1/p^\infty})= \Spec \F_p(t)[x_1,
x_2, \ldots]/ \langle x_1^p - t, x_i^p - x_{i-1} : i=2,3,\ldots
\rangle .$$ Since $\Omega_{Y/X}$ is generated as a
$\F_p(t^{1/p^\infty})$-vector space by $\{ dx_i : i=1,2,\ldots \}$ and
since the relation $x_{i+1}^p - x_i$ implies that $dx_i$ is zero, it
follows that $\Omega_{Y/X}=0$. Also, $Y \rightarrow X$ is clearly
profinite and flat. Since the field extension
$\F_p(t^{1/p^\infty})/\F_p(t)$ is purely inseparable, and since any
finite \'etale $X$-scheme is a finite disjoint union of spectra of
finite separable extensions of $\F_p(t)$, $Y$ is not an inverse limit
of finite \'etale $X$-schemes.

\begin{figure}[ht]
\begin{center}
\setlength{\unitlength}{0.00083333in}
\begingroup\makeatletter\ifx\SetFigFont\undefined%
\gdef\SetFigFont#1#2#3#4#5{%
  \reset@font\fontsize{#1}{#2pt}%
  \fontfamily{#3}\fontseries{#4}\fontshape{#5}%
  \selectfont}%
\fi\endgroup%
{\renewcommand{\dashlinestretch}{30}
\begin{picture}(1216,495)(0,-10)
\put(608.000,759.000){\arc{1500.000}{0.6435}{2.4981}}
\put(608.000,-441.000){\arc{1500.000}{3.7851}{5.6397}}
\put(983,384){\makebox(0,0)[lb]{{\SetFigFont{8}{9.6}{\rmdefault}{\mddefault}{\updefault}$q$}}}
\put(83,384){\makebox(0,0)[lb]{{\SetFigFont{8}{9.6}{\rmdefault}{\mddefault}{\updefault}$p$}}}
\end{picture}
}
\end{center}
\caption{An example showing that the notion
of profinite-\'etale is not Zariski-local}
\label{banana}
\end{figure}
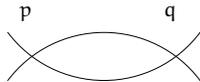

(b) Unlike covering spaces in topology, the property of being
profinite-\'etale is not Zariski-local on the target.  Here is an
example.  Consider the arithmetic genus $1$ complex curve $C$ obtained
by gluing two $\proj^1$'s together along two points, and name the
nodes $p$ and $q$ (Figure~\ref{banana}).  Consider the
profinite-\'etale covering space $Y \rightarrow C-p$ given by $$\Spec
\oh_{C-p} [\dots, x_{-1}, x_0, x_1, \dots] / (x_i^2-1)$$ and the
profinite-\'etale covering space $Z \rightarrow C-q$ given by $$\Spec
\oh_{C-q} [\dots, y_{-1}, y_0, y_1, \dots] / (y_i^2-1).$$  Glue $Y$ to
$Z$ (over $C$) by identifying $x_i$ with $y_i$ on the ``upper
component'', and $x_i$ with $y_{i+1}$ on the ``lower component''.
Then $Y \cup Z \rightarrow C$ is not profinite-\'etale, as it does not factor non-trivially through any finite \'etale morphism. To see this, suppose that we had a $C$ map $Y \cup Z \rightarrow W$ with $W \rightarrow C$ finite \'etale. We have a functor taking a $C$ scheme to the topological fiber over some point of the associated complex analytic space. Call this the fiber functor. Note that the fiber of $Y \cup Z$ is homeomorphic to $\F_2^{\Z}$ with the product topology. Let $s$ be the ``shift
indexing by $1$'' operator on the fiber, so $s(g(i)) = g(i+1)$ for $g: \Z \rightarrow \F_2$. Since finite \'etale covers are topological covering spaces, the fiber of $W$ has an action of $\Z$, and therefore must be of the form $\coprod_I \Z/m_i$ for some finite set $I$, positive integers $m_i$, and with $\Z$ acting by $+1$ on each $\Z/m_i$. Applying the fiber functor to $Y \cup Z \rightarrow W$ gives a $\Z$ equivariant continuous map $h: \F_2^{\Z} \rightarrow \coprod_I \Z/m_i$, where $\Z$ acts on $\F_2^{\Z}$ via $s$. We may assume that $h$ is surjective and show that $\sum m_i =1$, which is equivalent to showing $Y \cup Z \rightarrow W$ is a trivial factorization: if $\sum m_i > 1$, then $\F_2^{\Z}$ admits two disjoint, non-empty, open, $m \Z$-equivariant sets $U_1$ and $U_2$ for $m = \prod m_i$. For any $l,u$ in $\Z$, we have the map $$r_{(l,u)}: \F_2^{\Z} \rightarrow \F_2^{\{l,l+1,\ldots, u-1, u \}}$$ given by restriction of functions. By definition of the product topology, $U_i$ contains a subset $V_i$ of the form $V_i = r^{-1}_{(l_i,u_i)} (S_i) $ for some integers $l_i,u_i$ and a non-empty subset $S_i$ of $F_2^{\{l_i,l_i+1,\ldots, u_i-1, u_i \}}$. It follows that there is $N$ such that $s^{mN} V_1 \cap V_2 \neq \emptyset$, giving a contradiction.


A map from a connected $X$-scheme to a profinite-\'etale covering space of $X$ is determined by the image of a geometric point:

\begin{proposition}
{Let $(X,x)$ be a connected, geometrically-pointed scheme, and let $\wp:(Y,y) \rightarrow (X,x)$ be a profinite-\'etale covering space. If $f:(Z,z) \rightarrow (X,x)$ is a morphism from a connected scheme $Z$ and $\tilde{f}_1$ and $\tilde{f}_2$ are two lifts of $f$ taking the geometric point $z$ to $y$, then $\tilde{f}_1=\tilde{f}_2$. (A lift of $f$ means a map $\tilde{f}: Z \rightarrow Y$ such that $ \wp \circ \tilde{f}= f$.)
}\label{lifts_from_connected_determined_by_image_point_put_further_back}
\end{proposition}

\begin{proof}
By the universal property of the inverse limit, we reduce to the case where $\wp$ is finite \'etale. Since the diagonal of a finite \'etale morphism is an open and closed immersion, the proposition follows.
\end{proof}

Geometric points of a connected scheme lift to a profinite-\'etale covering space:

\begin{proposition}
{Let $X$ be a connected scheme, $x$ a geometric point of $X$, and $f:
  Y \rightarrow X$ a profinite-\'etale covering space. Then there is a
  geometric point of $Y$ mapping to $x$.
}\label{geom_pt_lift_through_pfe} \end{proposition}

\begin{proof}
Since $Y \rightarrow X$ is profinite-\'etale, $Y = \varprojlim_I Y_i$, where $Y_i \rightarrow X$ is finite \'etale and $I$ is a directed set. Let $\cF_x(Y_i)$ denote the geometric points of $Y_i$ mapping to $x$. Since $X$ is connected, $\cF_x(Y_i)$ is non-empty (because finite \'etale maps are open and closed \cite[\S 2 Th.\ 2.4.6]{ega4ii} \cite[\S 6 Prop.\ 6.1.10]{ega2} and induce finite separable extensions of residue fields of points \cite[\S 17 Th.\ 17.4.1]{ega4iv}). Since $Y_i \rightarrow X$ is finite, $\cF_x(Y_i)$ is finite. The set of geometric points of $Y$ mapping to $x$ is $\varprojlim_I \cF_x(Y_i)$. $\varprojlim_I \cF_x(Y_i)$ is non-empty because an inverse limit over a directed set of non-empty finite sets is non-empty \cite[Prop.\ I.I.4]{RZ}.
\end{proof}

\subsection{Example:  profinite sets give Yoneda trivial profinite-\'etale covering spaces}\label{trivial_pfe_scheme_associated_to_profinite_set}
If $S$ is a profinite set, define the {\em trivial $S$-bundle}
over $X$ by 
$$\underline{S}_X := \underline{\Spec} \;  \left( \Maps_{cts}(S, \oh_X)  \right)$$ where $\oh_X(U)$ is given the discrete topology for all open $U
\subset X$.  It is straightforward to verify that $\underline{S}_X \rightarrow X$ is a Yoneda trivial covering space with distinguished sections canonically homeomorphic to $S$, and that if $S = \varprojlim_I S_i$, then $\underline{S} = \varprojlim_I \underline{S_i}$. We will see that this example describes all Yoneda trivial profinite-\'etale covering spaces (Proposition \ref{Yoneda_triv_means_triv_bundle}).

The topology on the distinguished sections of a Yoneda trivial profinite-\'etale covering space is profinite:

\begin{proposition}{Let $f: Y \rightarrow X$ be a Yoneda trivial
    profinite-\'etale covering space with distinguished set of
    sections $S$. Let $p$ be any point of $\fT(X).$ Let
    $\cF_p(\fT(f))$ be the fiber of $\fT(f): \fT(Y) \rightarrow
    \fT(X)$ above $p$. The continuous map $S \rightarrow
    \cF_p(\fT(f))$ given by evaluation at $p$ is a homeomorphism. In
    particular, $S$ is profinite. }\label{top_S_is_top_fiber_hence_pf}
\end{proposition}

\begin{proof}
Since $f$ is profinite-\'etale, we may write $f$ as $\varprojlim f_i$ where $f_i: Y_i \rightarrow X$ is a finite \'etale covering space indexed by a set $I$. By \cite[\S 8 Prop.\ 8.2.9]{ega4iii}, the natural map $\fT(Y) \rightarrow \varprojlim \fT(Y_i)$ is a homeomorphism. Since $f_i$ is finite, $\cF_p(\fT(f_i))$ is finite. Thus, $\cF_p(\fT(f))$ is profinite.

For any $p' \in \cF_p(\fT(f))$, the extension of residue fields $k(p) \subset k(p')$ is trivial since the map $\Spec k(p') \rightarrow Y$ must factor through $X$ by Yoneda triviality. It follows that we have a unique lift of $\Spec k(p) \rightarrow X$ through $f$ with image $p'$. By definition of Yoneda triviality, we have that $p'$ is in the image of a unique element of $S$. Thus $S \rightarrow \cF_p(\fT(f))$ is bijective. $S \rightarrow \cF_p(\fT(f))$ is continuous, because $S$ is topologized by pointwise convergence.

Since $S$ is given the topology of pointwise convergence, to show that the inverse $\cF_p(\fT(f)) \rightarrow S$ is continuous is equivalent to showing that for any $q$ in $\fT(X)$, the map $$\cF_p(\fT(f)) \rightarrow S \rightarrow \cF_q(\fT(f))$$ is continuous.

The set of sections $S$ produces a set sections $S_i$ of $f_i$. Since $Y \rightarrow Y_i$ is profinite-\'etale, $Y \rightarrow Y_i$ is integral. Thus,  $\cF_p(\fT(f)) \rightarrow \cF_p(\fT(f_i))$ is surjective. It follows that for any $p_i' \in \cF_p(\fT(f_i))$, $p_i'$ is in the image of one of the sections in $S_i$ and that $k(p_i')=k(p)$. Since $Y_i \rightarrow X$ is
finite-\'etale and $X$ is connected, it follows that $Y_i \cong
\coprod_{S_i} X.$ The isomorphisms $Y_i \cong
\coprod_{S_i} X$ identify $\cF_p(\fT(f))$, $S$, and $\cF_q(\fT(f))$ with $\varprojlim S_i$ compatibly with the evaluation maps.  
\end{proof}

Yoneda trivial profinite-\'etale covering spaces are trivial $S$-bundles, where $S$ is the distinguished set of sections as a topological space. In fact, taking such a covering space to its distinguished sections is an equivalence of categories:

\begin{proposition}{Let $X$ be a connected scheme and let $f: Y \rightarrow X$ be a Yoneda trivial profinite-\'etale covering space. Let $S$ denote the distinguished set of sections of $f$. Then there is a canonical isomorphism of $X$-schemes $Y \cong \underline{S}_X$.
Furthermore, if $f_1: Y_1 \rightarrow X$ and $f_2: Y_2 \rightarrow X$ are two Yoneda trivial profinite-\'etale covering spaces with distinguished sets of sections $S_1$ and $S_2$ respectively, then the map $$\Maps_{cts}(S_1,S_2) \rightarrow \Maps_X(Y_1,Y_2)$$ induced by $\Maps_{cts}(S_2, \oh_X) \rightarrow \Maps_{cts}(S_1, \oh_X) $ is a bijection. ($\Maps_X(Y_1,Y_2)$ denotes the set of scheme morphisms $Y_1 \rightarrow Y_2$ over $X$.)
 }\label{Yoneda_triv_means_triv_bundle}\end{proposition}

\begin{proof}
Since every element of $S$ is a map $X \rightarrow Y$, we have a canonical map $ S \times \mathcal{O}_Y \rightarrow \mathcal{O}_X$. 
(By $S \times \mathcal{O}_Y$, we mean a product of copies of $\mathcal{O}_Y$ indexed by $S$.)
By adjointess, we have  $\mathcal{O}_Y \rightarrow \Maps(S, \mathcal{O}_X)$.

Since $f$ is profinite-\'etale, there is an inverse system of finite
\'etale $X$-schemes $\{ Y_i \rightarrow X \}_{i \in I}$ such that $Y
\cong \varprojlim_I Y_i$. As in the proof of Proposition \ref{top_S_is_top_fiber_hence_pf}, for each $i \in I$, $S$ induces a (finite) set of sections $S_i$ of $Y_i \rightarrow X$ and, furthermore, $Y_i \cong \coprod_{S_i} X$ and $S \cong \varprojlim_I S_i.$

Since $Y \cong \varprojlim_I Y_i$, the
map $\mathcal{O}_Y \rightarrow \varinjlim_I \Maps(S_i,
\mathcal{O}_X)$ is an isomorphism. Note that $\varinjlim_I
\Maps(S_i, \mathcal{O}_X) = \Maps_{\rm{cts}}(\varprojlim_I
S_i, \mathcal{O}_X)$. Thus we have a canonical isomorphism of $X$-schemes $Y = \underline{S}_X$.

Now consider $f_1$ and $f_2$. Given $g \in  \Maps(Y_1,Y_2)$ and $s_1 \in S_1$, we have a section $g \circ s_1$ of $f_2$, and therefore an element $s_2 \in S_2$. Thus $g$ determines a map $S_1 \rightarrow S_2$. Since  the evaluation maps $S_j \rightarrow \cF_p(\fT(f_j))$ $j=1,2$ and the map $\fT(g): \cF_p(\fT(f_1)) \rightarrow \cF_p(\fT(f_2))$ fit into the commutative diagram

$$
 \xymatrix{ S_1 \ar[r] \ar[d] & S_2 \ar[d]\\
   \cF_p(\fT(f_1)) \ar[r] &  \cF_p(\fT(f_2)),}
$$

\noindent the map $S_1 \rightarrow S_2$ is continuous by Proposition \ref{top_S_is_top_fiber_hence_pf}. We therefore have $\Maps(Y_1,Y_2) \rightarrow  \Maps_{cts}(S_1,S_2)$.

For $s$ in $S$, the map $s: X \rightarrow Y_2$ is identified with the map $X \rightarrow \underline{S}_X$ induced by the `evaluation at $s$' map $ \Maps_{cts}(S, \mathcal{O}_X) \rightarrow \mathcal{O}_X$, under the isomorphism $Y = \underline{S}_X$. It follows that \begin{equation*}\Maps_{cts}(S_1,S_2) \rightarrow \Maps(Y_1,Y_2)
\rightarrow \Maps_{cts}(S_1,S_2)\end{equation*} is the identity.  Likewise, for $s_1$ in $S_1$, the composition $$ \xymatrix{\Maps(Y_1,Y_2) \ar[r]&
\Maps_{cts}(S_1,S_2) \ar[r] & \Maps(Y_1,Y_2) \ar[r]^{s_1^*} &  \Maps(X,Y_2)}$$ is given by $$g \mapsto g \circ s_1. $$

Because
$\coprod_{s_1 \in S_1} s_1: \coprod_{S_1} X \rightarrow Y_1$ is an
fpqc cover, $$ \prod_{s_1 \in S_1} s_1^*:  \Maps(Y_1,Y_2) \rightarrow \prod_{S_1} \Maps( X, Y_2) $$ is injective, and it follows that $$\Maps(Y_1,Y_2) \rightarrow
\Maps_{cts}(S_1,S_2) \rightarrow \Maps(Y_1,Y_2)$$ is the
identity. \end{proof}

Heuristically, an object is Galois if it has maximal symmetry. Since automorphisms $\Aut(Y/X)$ of a covering space $Y \rightarrow X$ are sections of the pullback $Y \times_X Y \rightarrow Y,$ it is reasonable to define a covering space to be Galois if the pullback is Yoneda trivial: 

\begin{definition}\label{Galoisdef} A profinite-\'etale covering space $Y \rightarrow X$ is defined to be {\em Galois} if the left (or equivalently, right) projection $Y \times_X Y \rightarrow Y$ is Yoneda trivial. \end{definition}

For a Galois covering space $Y \rightarrow X$ with $Y$ connected, $\Aut(Y/X)$ is a profinite group; the topology on $\Aut(Y/X)$ comes from identifying $\Aut(Y/X)$ with the space of distinguished sections of $Y \times_X Y \rightarrow Y$ and applying Proposition \ref{top_S_is_top_fiber_hence_pf}. 

\subsection{Example: Trivial profinite group schemes over $X$}
\label{trivialprofinitegroup} 
If $G$ is a
profinite group with inverse $i$ and multiplication $m$, define the
{\em trivial $G$-bundle} as the $X$-scheme $\underline{G}_X$ of Example \ref{trivial_pfe_scheme_associated_to_profinite_set} with the
following group scheme structure.
We describe a Hopf algebra structure over an open set $U$;
this construction will clearly glue to yield a sheaf of Hopf algebras.
The coinverse map sends 
\begin{equation}
\label{star}\xymatrix{G \ar[r]^-f_-{\rm{cts}} & \oh_X(U)}
\end{equation}
to the composition
$$
\xymatrix{G \ar[r]^{i} & G \ar[r]^-f_-{\rm{cts}} &  \oh_X(U).}
$$
The coinverse $f \circ i$ is indeed continuous, as it is the composition of two
continuous maps.  The comultiplication map sends \eqref{star} to the
composition
\begin{equation}
\label{comult}
\xymatrix{ G \times G 
\ar@/^/[rr]
\ar[r]_-{m} & G \ar[r]_-{f \; \; \rm{cts}} &  \oh_X(U)}
\end{equation}
using the isomorphism $$\Maps_{\rm{cts}}(G \times G, \oh_X(U)) \cong
\Maps_{\rm{cts}}(G, \oh_X(U)) \otimes_{\oh_X(U)} \Maps_{\rm{cts}}(G,
\oh_X(U))$$ where $G \times G$ has the product topology.
The map \eqref{comult} is
continuous, as it is the composition of two continuous maps.  The
coidentity map is the canonical map
$\Maps_{\rm{cts}}(G, \oh_X) \rightarrow \oh_X$ given by
evaluation at $e$.  The fact that these maps satisfy the axioms of a
Hopf algebra is the fact that $(G, e, i, m)$ satisfies the
axiom of a group.

The trivial $G$-bundle on any $X$ is clearly pulled back from 
the trivial $G$-bundle on $\Spec \Z$.  

\subsection{Example:  $\hat{\Z}$, roots of unity, and Cartier duality}The following example is well-known.  It is included because it is an explicit example of the construction of \S \ref{trivialprofinitegroup}.  

The roots of unity form a Hopf algebra: let $A$ be a ring and define
$$
A[\mu_{\infty}] := A [t_{1!}, t_{2!}, t_{3!}, \dots] / (t_{1!}-1,
t^2_{2!} -t_1, \dots, t^n_{n!} - t_{(n-1)!}, \dots).$$  Give
$A[\mu_{\infty}]$  a Hopf algebra structure  by coinverse  
$\iota:  A[\mu_{\infty}] \rightarrow A[\mu_{\infty}]$ given by 
$\iota: t_n
\mapsto t_n^{-1}$ and comultiplication $\mu:
  A[\mu_{\infty}] \rightarrow A[\mu_{\infty}]\otimes_A
  A[\mu_{\infty}]$ given by  
$\mu: t_n \mapsto t_n' t_n''$.  (The new
variables $t'_n$ and $t_n''$ are introduced to try to avoid notational
confusion:
 they are new names
for the $t_n$-coordinates on  the left and right factors respectively of
$A[\mu_{\infty}]\otimes_A A[\mu_{\infty}]$.)\label{bubblegum}

Let $A$ be a ring containing a primitive $n^{th}$ root of unity for
any positive integer $n$. (In particular $\chr A = 0$.)  The $t_{j!}$
correspond to continuous characters $\hat{\Z} \rightarrow A^*.$ For
example, $t_2$ corresponds to the continuous map sending even elements
to $1$ and odd elements to $-1$ (i.e.\ $n \mapsto (-1)^n$). (Choosing
such a correspondence is equivalent to choosing an isomorphism between
$\hat{\Z}$ and $\mu_{\infty}(A)$.) The hypothesis $\chr A =0$ implies
that $A [\mu_{\infty}]$ is isomorphic to the subalgebra of continuous
functions $\hat{\Z} \rightarrow A $ generated by the continuous
characters. Because the characters span the functions $\Z/n
\rightarrow A$, it follows that $$ \underline{\hat{\Z}}_{\Spec A}
\cong \Spec A [\mu_{\infty}] .$$ Such an isomorphism should be
interpreted as an isomorphism between $\underline{\hat{\Z}}$ and its
Cartier dual.  \label{Zhat}

Combining Proposition \ref{Yoneda_triv_means_triv_bundle} and Example \ref{trivialprofinitegroup} shows that a connected Galois covering space pulled back by itself is the trivial group scheme on the automorphisms:

\begin{proposition}
{Let $f:Y \rightarrow X$ be a Galois profinite-\'etale covering space with $Y$ connected. Then 
 \begin{equation}
\label{Galois_fpqc_cover}
\xymatrix{ \underline{\Aut(Y/X)}_X \times_X Y \ar[rr]^-{\mu} \ar[d] && Y \ar[d] \\
Y \ar[rr] && X}
\end{equation}
\noindent is a fiber square such that the map $\mu$ is an
action.}\label{trivial_Aut_bundle_Galois_covering_fiber_square} \end{proposition}

\section{Algebraic universal covers}

\label{s:UC}

\begin{definition}\label{sc_def} A connected scheme $X$ is {\em simply
  connected} if all covering spaces are Yoneda trivial.  With covering space defined as in Definition \ref{d:covering_space}, this is equivalent to the usual definition that $X$ is simply connected if a connected finite \'etale $X$-scheme is isomorphic to $X$ (via the structure map).\end{definition}

\subsection{Remark}  Of course, many other similarly parsimonious definitions are possible, so we gives some indication of the advantages of this one.   As with many other ``functor of points''
style definitions in algebraic geometry, this particular definition makes constructions and proof simpler.  One could define $X$ to be simply connected if any connected cover of $X$ is an isomorphism.   But any definition of $X$ being simply connected which only places restrictions on connected covers of $X$ will run into the difficulty that when $X$ is not locally Noetherian, one can not always express a cover as a disjoint union of connected components.

\begin{definition}\label{UC_def} A covering space $p: \tX \rightarrow
X$ of a connected scheme $X$ is a {\em universal cover} if $\tX$ is
connected and simply connected.\end{definition}

\begin{proposition} {Suppose $p: \tX \rightarrow X$ is a
  universal cover of a connected scheme $X$, $\tilde{x}$ is a geometric point
  of $\tX$, and $x = p(\tilde{x})$ is the corresponding geometric
  point of $X$.  Then $(\tX, \tilde{x}) \rightarrow (X, x)$ is initial
  in the category of geometrically pointed covering spaces of
  $(X,x)$.}\label{newbie}\end{proposition}

\begin{proof}
Suppose $(Y,y) \rightarrow (X,x)$ is a geometrically pointed covering
space of $(X,x)$.   Since covering spaces are stable under pullback,
$\tX \times_X Y \rightarrow \tX$ is a profinite-\'etale covering
space.  Since $\tX$ is simply connected, 
$\tX \times_X Y \rightarrow \tX$ is Yoneda trivial.  In particular, $\tX \times_X Y \rightarrow \tX$
admits a section sending $\tilde{x}$ to $\tilde{x} \times y$, from
which we have a map  of covering spaces $(\tX, \tilde{x}) \rightarrow (Y,y)$, which is
unique by the connectedness of $\tX$. \end{proof}

\begin{proposition}
{Let $X$ be a connected scheme. Then a universal cover of $X$ is unique up to (not necessarily unique) isomorphism.}\label{UC_unique_up_to_non_unique_iso}\end{proposition}

\begin{proof}
Let $\tX_1,\tX_2$ be two universal covers of $X$. Since covering
spaces are stable under pull-back, $\tX_1 \times_X \tX_2 \rightarrow
\tX_1$ is a profinite-\'etale covering space. Since $\tX_1$ is simply
connected, $\tX_1 \times_X \tX_2 \rightarrow \tX_1$ is Yoneda trivial
and thus admits a section, whence we have a map of $X$-schemes $f: \tX_1 \rightarrow \tX_2.$ We see that $f$ is an isomorphism as follows: since $\tX_1 \times_X \tX_2 \rightarrow \tX_2$ is Yoneda trivial, the map $\textrm{id} \times f: \tX_1 \rightarrow \tX_1 \times_X \tX_2$ factors through $f: \tX_1 \rightarrow \tX_2$ by a distinguished section $g \times \textrm{id}$ of $\tX_1 \times_X \tX_2 \rightarrow \tX_2.$ In particular $gf: \tX_1 \rightarrow \tX_1$ is the identity. Switching the roles of $\tX_1$ and $\tX_2$ we can find $f': \tX_1 \rightarrow \tX_2$ such that $f'g:\tX_2 \rightarrow \tX_2$ is the identity. Thus $f = f' g f = f',$ and we have that $f$ is an isomorphism with inverse $g$.
\end{proof}

\begin{proposition}
{Let $X$ be a connected qcqs scheme equipped with a geometric point $x$. Suppose $p:(\tX,\tilde{x}) \rightarrow (X,x)$  is an initial object among geometrically-pointed covering spaces of $X$ such that $\tX$ is connected. Then $\tX$ is
a simply connected
 Galois covering space, and in particular is a universal cover.  Thus
 universal covers are Galois.}\label{initial_pointed_covering_space}\end{proposition}

Somewhat conversely, a connected, simply connected, Galois covering space is an
initial object among connected geometrically-pointed covering spaces
of $X$, by Proposition~\ref{newbie}.

\begin{proof} We first show that $\tX$ is simply connected.  Let $q: \tY
\rightarrow \tX$ be a covering space of $\tX$ and let $S$ be the set
of sections of $q$. We will show that $q$ is Yoneda trivial with
distinguished set of sections $S$. Let $Z$ be a connected
$\tX$-scheme. We need to show that $S \rightarrow \Maps_{\tX}(Z,\tY)$
is bijective. Injectivity follows from Proposition
\ref{lifts_from_connected_determined_by_image_point_put_further_back}. From
Proposition
\ref{lifts_from_connected_determined_by_image_point_put_further_back}
it also follows that we may assume that $Z \rightarrow \tX$ is a
geometric point of $\tX.$ Let $z$ be a geometric point of $\tX$. By Proposition \ref{geom_pt_lift_through_pfe}, we may lift $z$ to a geometric point $\tilde{z}$ of $\tY.$ Applying Proposition
\ref{lifts_from_connected_determined_by_image_point_put_further_back}
again, we see that it suffices to construct a map of $X$-schemes
$(\tX, z) \rightarrow (\tY,\tilde{z})$. Since $X$ is qcqs,
profinite-\'etale maps are closed under composition. Thus $\tY$ is an inverse limit of finite
\'etale $X$-schemes. Thus by Proposition
\ref{lifts_from_connected_determined_by_image_point_put_further_back},
it suffices to show that for any pointed finite \'etale $(Y,y)
\rightarrow (X,pz)$, we have an $X$ map $(\tX, z) \rightarrow
(Y,y)$. Take $Y \rightarrow X$ finite \'etale, and let $d$ be the
degree of $Y$. Since $p$ is initial, we have $d$ maps $\tX \rightarrow
Y$ over $X$. By Proposition
\ref{lifts_from_connected_determined_by_image_point_put_further_back},
we therefore have an $X$ map $(\tX, z) \rightarrow (Y,y)$. Thus $\tX$
is simply connected.

Since $\tX$ is simply connected, $\tX \times_X \tX \rightarrow \tX$ is Yoneda trivial, and therefore $\tX$ is a Galois covering space. \end{proof}

\begin{proposition} {If $X$ is a connected qcqs scheme, then a
  universal cover $p: \tX \rightarrow X$ exists.} \label{UCexists}
\end{proposition}

\subsection{Remark on Noetherian conditions} If $X$ is Noetherian, in
general $\tX$ will not be Noetherian. We will see (Theorem~\ref{etalefglink}) that the
geometric fibers of $p$ are in natural bijection with the \'etale
fundamental group.  Thus if $X$ has infinite \'etale fundamental
group, and a closed point $q$ with algebraically closed residue field, then
$p^{-1}(q)$ is dimension $0$ (as $p$ is integral) with an infinite
number of points, so $\tX$  has a closed subscheme which is not Noetherian
and is thus  not Noetherian itself.  However, such a solenoid is not so
pathological. For example, the \'etale topological type as in \cite{am} presumably extends to universal covers of locally Noetherian schemes, although we have not worked this out carefully. Also, by \cite[Pt.~0, Lem~ 10.3.1.3]{ega3}, the local rings of $\tX$
are Noetherian.
\label{Nremark}

\begin{proof} {\em (of Proposition \ref{UCexists})}
Choose a geometric point $x: \Spec \Omega \rightarrow X$. By Proposition \ref{initial_pointed_covering_space}, it suffices to show that the category of pointed covering spaces of $(X,x)$ has a connected initial object.

If $(Y_\nu, y_\nu)$ are two geometrically-pointed connected finite \'etale
$(X,x)$-schemes, we will say that $(Y_2, y_2) \geq (Y_1, y_1)$ if
there is a morphism of pointed $(X,x)$-schemes $(Y_2, y_2) \rightarrow (Y_1,
y_1)$.
The diagonal of a finite \'etale map is an open and closed
immersion, so an $X$-map from a
connected scheme to a finite \'etale $X$-scheme is determined by the
image of a single geometric point. Thus the symbol $\geq$ is a
partial order on isomorphism classes of connected pointed finite
\'etale $X$-schemes.

It is straightforward to see that the isomorphism classes of connected finite etale $X$-schemes form a set, for any scheme $X$.  Indeed, the isomorphism classes of affine finite type $X$-schemes (i.e.\ schemes with an affine finite type morphism to $X$)
 form a set.
First show this for each affine scheme.  (Each affine finite  type  $X$-scheme
can be described in terms of a fixed countable list of variables.  The possible
relations form a set, and the actual relations lie in the power set of this set.)
For each pair of affine opens $U_i$ and $U_j$, cover 
$U_i \cap U_j$ with a set of affine opens $U_{ijk}$ simultaneously distinguished in $U_i$ and $U_j$.
For each $ijk$, and each affine   finite type cover of $U_i$ and $U_j$, there is a set of morphisms from the restriction
of the cover of $U_i$ to the restriction of the cover of $U_j$ (look at the corresponding rings, and choose images of generators).
Within this set of data, we take the subset where these morphisms (for each $ijk$) are isomorphisms; then take the subset where these morphisms
glue together (yielding a affine  finite type cover of $X$).  Then quotient by isomorphism.

The set $I$ of isomorphism classes of connected finite \'etale
$X$-schemes equipped with $\geq$ is directed: suppose $(Y_1, y_1)$ and
$(Y_2, y_2)$ are two geometrically-pointed connected $(X,x)$-schemes.
Then $(Y_1 \times_X Y_2, w := y_1 \times y_2)$ is a
geometrically-pointed finite \'etale $(X,x)$-scheme.  Even though we have
made no Noetherian assumptions, we can make sense of ``the connected
component $Y'$ of $Y_1 \times Y_2$ containing $w$''.  If $Z
\rightarrow X$ is a finite \'etale cover, then it has a well-defined
degree, as $X$ is connected.  If $Z$ is not connected, say $Z = Z_1
\coprod Z_2$, then as $Z_i \rightarrow X$ is also finite \'etale
($Z_i$ is open in $Z$ hence \'etale over $X$, and closed in $Z$, hence
finite), and has strictly smaller degree.  Thus there is a smallest
degree $d$ such that there exists an open and closed $W
\hookrightarrow Y_1 \times_X Y_2$ containing $y_1 \times y_2$ of
degree $d$ over $X$, and $W$ is connected.  
(Note that $W$ is unique:  the set of such $W$ is closed
under finite intersections, and the intersection of two such, say $W_1$
and $W_2$, has degree strictly less than that of $W_1$ and $W_2$.)  Then $(W,w) \geq (Y_i,
y_i)$.

By \cite[\S 8 Prop.\ 8.2.3]{ega4iii}, inverse limits with affine
transition maps exist in the category of schemes, and the inverse
limit is the affine map associated to the direct limit of the sheaves
of algebras. Define $\tilde{X} :=
\varprojlim_I Y_i$, where we have chosen a representative pointed
connected finite \'etale $X$-scheme $(Y_i, y_i)$ for each $i\in
I$. The geometric points $\{ y_i \}_{i\in I}$ give a canonical
geometric point $\tilde{x}$ of $\tilde{X}$.

By \cite[\S 8 Prop.\ 8.4.1(ii)]{ega4iii}, since $X$ is quasicompact,
$\tilde{X}$ is connected. (This is the only place in the proof that the category of pointed covering spaces of $(X,x)$ has a connected initial object where the
  quasicompactness hypotheses is used.)
  
$(\tilde{X}, \tilde{x})$ admits a map to any pointed finite \'etale $(X,x)$-scheme by construction. This map is unique because $\tilde{X}$ is connected. Passing to the inverse limit, we see that $(\tilde{X}, \tilde{x})$ is an initial object in pointed profinite-\'etale $X$-schemes. 
\end{proof}

\begin{corollary}{Profinite-\'etale covering spaces of connected qcqs schemes are profinite-\'etale locally (i.e.\ after pullback to a profinite-\'etale cover) Yoneda trivial.}\label{pfe_cov_spaces_pfe_loc_Yon_triv}\end{corollary}
 
The remainder of this section is devoted to examples and properties of universal covers.  It is not necessary for the construction of the fundamental group family of \S \ref{s:FG}. 

\subsection{Universal covers of group schemes}\label{UCgrp_sch_b_point}
The following result and proof are the same
as for Lie groups.

\begin{theorem} {Let $X$ be a connected group scheme finite type over an
  algebraically closed field $k$.  Suppose $\chr k = 0$
  or $X$ is proper.  Choose any preimage $\tilde{e} \in \tX$ above $e
  \in X$.  Then there exists a unique group scheme structure on $\tX$
  such that $\tilde{e}$ is the identity and $p$ is a morphism of group
  schemes over $k$.} \label{UCgroupscheme} \end{theorem}
 
The choice of $\tilde{e}$ is not important: if $\tilde{e}'$ is another
choice, then $(\tilde{X}, \tilde{e}) \cong (\tilde{X}, \tilde{e}')$.  If $k$ is not algebraically closed and $\chr k =0$, then $\tX$ is the universal cover of $X_{\kbar}$, and we can apply Theorem~\ref{UCgroupscheme} to $X_{\kbar}$, obtaining a similar statement, with a more awkward wording.  For example, the residue field of
$\tilde{e}$ is the algebraic closure of that of $e$. To prove
Theorem~\ref{UCgroupscheme}, we use a well-known lemma whose proof we include due to lack of a reference.

\begin{lemma} {Suppose $X$ and $Y$ are connected finite type
  schemes over an algebraically closed field $k$.  Suppose $\chr k =
  0$ or $X$ is proper.  Then $\tilde{X} \times \tilde{Y}$ is simply
  connected.  Equivalently, a product of universal covers is naturally
  a universal cover of the product.}\label{janlemma}\end{lemma}

\begin{proof} This is equivalent to the following statement about the \'etale
fundamental group.  Suppose $X$ and $Y$ are finite type over an
algebraically closed field $k$, with $k$-valued
points $x$ and $y$ respectively. 
Suppose $X$ is proper or $\chr k =0$.  Then the natural group homomorphism
$$
\pi_1^{et}(X \times Y, x \times y) \rightarrow \pi_1^{et} (X,x) \times \pi_1^{et}(Y,y)
$$
is an isomorphism. The characteristic $0$ case follows by reducing to $k = \C$ using the Lefschetz principle, and comparing $\pi_1^{et}$ to the topological fundamental group. The $X$ proper case is \cite[Exp.\ X Cor.\
1.7]{sga1}.
\end{proof}

\begin{proof} {\em (of Theorem \ref{UCgroupscheme})} We first note the following: suppose $(W,w) \rightarrow (Y,y)$ is a geometrically pointed covering space. If we have a map of geometrically pointed schemes $f:(Z,z) \rightarrow (Y,y)$ from a simply connected scheme $Z$, then there is a unique lift of $f$ to a pointed morphism $\tilde{f}: (Z,z) \rightarrow (W,w),$ because $W \times_Y Z \rightarrow Z$ is a Yoneda trivial covering space.

Thus, there is a unique lift
  $\tilde{i}: \tX \rightarrow \tX$ lifting the inverse map $i: X
  \rightarrow X$ with $\tilde{i}(\tilde{e}) = \tilde{e}$.  By Lemma~\ref{janlemma}, $\tX \times \tX$ is simply
connected. Thus, there is a unique lift $\tilde{m}: \tX \times \tX \rightarrow \tX$
  of the multiplication map $m: X \times X \rightarrow X$ with
  $\tilde{m}( \tilde{e},\tilde{e}) = \tilde{e}$. It is straight forward to check that
  $(\tX, \tilde{e}, \tilde{i}, \tilde{m})$ satisfy the axioms of a
  group scheme.  For instance, associativity can be verified as follows:  we
  must show that $\tX \times \tX \times\tX \rightarrow \tX$ given by
  $((ab)c) (a(bc))^{-1}$ is the same
  as the identity $\tilde{e}$. Since associativity holds for $(X, e, i,m)$, both of these maps lie above $e: X \times X \times X \rightarrow X$. Since both send $\tilde{e} \times \tilde{e} \times \tilde{e}$ to $\tilde{e},$ they are equal. \end{proof}
  
 The assumption that $\chr k = 0$ or $X$ is proper is necessary for Theorem \ref{UCgroupscheme}, as shown by the following example of David Harbater (for which
we thank him).
 
\begin{proposition} {Let $k$ be a field with $\chr k = p >0$, and assume that $p$ is odd. The group law on $\G_a$ over $k$ does not lift to a group law on the universal cover.}\label{tilde(G_a)_not_grp_chr_p}\end{proposition}

\begin{proof}
Since the universal covers of $\G_a$ over $k$ and $k^s$ are isomorphic, we may assume that $k$ is separably closed. Let $p: \tilde{\G}_a \rightarrow \G_a$ denote the universal cover. Let $\tilde{e}$ be a $k$ point of $\tilde{\G}_a$ lifting $0$. Let $a: \G_a \times \G_a \rightarrow \G_a$ denote the addition map. Since $\tilde{\G}_a$ is simply connected, the existence of a commutative diagram

$$
\xymatrix{ \tilde{\G}_a \times_k  \tilde{\G}_a \ar[d]^{p \times p} \ar[rr] && \tilde{\G}_a \ar[d]^{p}  \\
 \G_a \times_k \G_a \ar[rr]^-{a} && \G_a}
$$

\noindent would imply that the composite homomorphism $a_* (p \times p)_*: \pi_1^{et}( \tilde{\G}_a \times_k  \tilde{\G}_a) \rightarrow \pi_1^{et}( \G_a \times_k \G_a) \rightarrow \pi_1^{et}( \G_a )$ is constant. (Say the basepoints of $\pi_1^{et}( \tilde{\G}_a \times_k  \tilde{\G}_a)$, $ \pi_1^{et}( \G_a \times_k \G_a)$, and  $\pi_1^{et}( \G_a )$ are $\tilde{e} \times \tilde{e}$, $e \times e$, and $e$ respectively, although this won't be important.) Thus, it suffices to find an element of $\pi_1^{et}( \tilde{\G}_a \times_k  \tilde{\G}_a)$ with non-trivial image under $a_* (p \times p)_*$.

Adopt the notation $\G_a \times_k \G_a = \Spec k[x] \times_k \Spec k[y]$. Consider the finite \'etale cover $W$ of  $\G_a \times_k \G_a=  \Spec k[x,y]$ given by the ring extension $k[x,y] \rightarrow k[x,y,w]/\langle w^p - w -xy \rangle $. $W$ is not pulled back from a finite \'etale cover of $\G_a$ under either projection map, as one sees with the Artin-Schreier exact sequence. It then follows from degree considerations that  $W$ and $\tilde{\G}_a \times_k  \tilde{\G}_a$ are linearly disjoint over $\G_a \times_k \G_a$. (By this we mean that the fields of rational functions of $W$ and $\tilde{\G}_a \times_k  \tilde{\G}_a$ are linearly disjoint over the field of rational functions of $\G_a \times_k \G_a$.) Thus $$\Aut((\tilde{\G}_a \times_k  \tilde{\G}_a \times_k W) / (\G_a \times_k \G_a)) \cong \Aut((\tilde{\G}_a\times \tilde{\G}_a)/ (\G_a \times \G_a)) \times \Aut(W/\G_a).$$

 It follows that there exists an element $\gamma$ of $\pi_1^{et}( \tilde{\G}_a \times_k  \tilde{\G}_a)$ which acts non-trivially on the fiber of $W$ pulled back to $\tilde{\G}_a \times_k  \tilde{\G}_a$.
 
 Choosing a lift of $\tilde{e} \times \tilde{e}$ to the universal cover of $\G_a \times \G_a$ allows us to view $\gamma$ as acting on any finite \'etale ring extension of $k[x,y]$ (by pushing $\gamma$ forward to $\G_a \times_k \G_a$ and using the isomorphism between elements of $\pi_1^{et}$ and automorphisms of the universal cover corresponding to a lift of base point). $\gamma$ therefore determines an automorphism of $k[x,y,w]/\langle w^p - w -xy \rangle$, and by construction, $\gamma$ acts non-trivially on $w$. Similarly, $\gamma$ determines automorphisms of $k[x,y,w_1]/\langle w_1^p - w_1 -x^2 \rangle$ and $k[x,y,w_2]/\langle w_2^p - w_2 -y^2 \rangle$, and in this case, $\gamma$ acts trivially on $w_1$ and $w_2$. 
 
Lifting $w$, $w_1$, and $w_2$ to functions on the universal cover, we have the function $z = 2w + w_1 + w_2$. Since $p \neq 2$, $\gamma$ acts non-trivially on $z$. Because $z$ satisfies the equation $$z^p -z = (x+y)^2,$$ $\gamma$ must act non-trivially on the fiber of the cover corresponding to $k[x,y] \rightarrow k[x,y,z]/\langle z^p - z -(x+y)^2 \rangle $.  Let $Z$ denote this cover.

Since $Z$ can be expressed as a pull-back under $a$, we have that $a_* (p \times p)_*\gamma$ is non-trivial.
 \end{proof}

\subsection{Examples} The universal cover can be described explicitly in a number of cases.
Of course,  if $k$ is
a field, then $\Spec k^s \rightarrow \Spec k$ is a universal cover.

\subsection{$\G_m$ over a characteristic $0$ field $k$}
This construction is also well known.
The Riemann-Hurwitz formula implies that the finite
\'etale covers of $\Spec k[t,t^{-1}]$ are obtained by adjoining roots of $t$ and
by extending the base field $k$.
Thus a universal cover is
$$
p: \Spec \overline{k} [t^{\Q}] \rightarrow \Spec k[t^{\Z}].$$ The
group scheme structure on the universal cover
(Theorem~\ref{UCgroupscheme}) is described in terms of the Hopf
algebra structure on $\overline{k} [t^{\Q}]$ given by coinverse $\iota: t^q \mapsto
t^{-q}$ and comultiplication $\mu: t^q \mapsto (t')^q (t'')^q$ ($q \in \Q$), which clearly lifts the group
scheme structure on $\G_m$. 
(Cf.\ Example~\ref{bubblegum}; analogously to there, the new
variables $(t')^q$ and $(t'')^q$, are introduced to try to avoid notational
confusion:
 they are new names
for the coordinates on  the left and right factors respectively of
$\overline{k}[t^{\Q}] \otimes_{ k[t^{\Z}]} \overline{k}[t^{\Q}]$.
Thus for example $t' = t''$, but $(t')^{1/2} \neq (t'')^{1/2}$.)
 Note that the universal cover is not
Noetherian.

\subsection{Abelian varieties}
We now explicitly describe the universal cover of an abelian variety over a field $k$.
We begin with separably  closed  $k$ for simplicity.

If $X$ is proper over a separably closed $k$, by the main
theorem of \cite{pardini}, the connected (finite) Galois covers
with {\em abelian} Galois group $G$ correspond to inclusions $\chi:
G^{\vee} \hookrightarrow \Pic X$, where $G^{\vee}$ is the dual group
(noncanonically isomorphic to $G$). 
The cover
corresponding to $\chi$ is
$\underline{\Spec} \oplus_{g \in G^{\vee}} \cL_{\chi(g)}^{-1}$
where
$\cL_{\chi}$ is the invertible sheaf corresponding to $\chi \in \Pic X$.

If $A$ is an abelian variety over $k$, then all Galois 
covers are abelian.  
Thus
$$
\tilde{A} = 
\underline{\Spec} \oplus_{\chi \text{ torsion}} \cL_{\chi}^{-1}
$$
where the sum is over over the torsion elements of $\Pic X$.

By Theorem~\ref{UCgroupscheme}, $\tA$ has a unique group scheme
structure lifting that on $A$ once a lift of the identity is chosen, when $k$ is  algebraically closed. In fact, $\tA$ has this group scheme structure with $k$ separably closed, and we now describe this
explicitly. Let $i: A \rightarrow A$ and $m: A \times A \rightarrow A$
be the inverse and multiplication maps for $A$.  
Then the
inverse map $\tilde{i}: \tA \rightarrow \tA$ is given by
$$
\xymatrix{\tilde{i}: \underline{\Spec} \oplus_{\chi \text{ torsion}} \cL_{\chi}^{-1}
\ar[r] & \underline{\Spec} \oplus_{\chi \text{ torsion}} i^* \cL_{\chi}^{-1}}
$$

\noindent using the isomorphism $i^* \cL \cong \cL^{-1}$ (for torsion sheaves,
by the Theorem of the Square).  The multiplication map $\tilde{m}: \tA
\times \tA \rightarrow \tA$ is via $m^* \cL \cong \cL \boxtimes
\cL$ (from the Seesaw theorem).

\begin{figure}[ht]
\begin{center}
\setlength{\unitlength}{0.00083333in}
\begingroup\makeatletter\ifx\SetFigFont\undefined%
\gdef\SetFigFont#1#2#3#4#5{%
  \reset@font\fontsize{#1}{#2pt}%
  \fontfamily{#3}\fontseries{#4}\fontshape{#5}%
  \selectfont}%
\fi\endgroup%
{\renewcommand{\dashlinestretch}{30}
\begin{picture}(2550,1839)(0,-10)
\put(1875,1362){\makebox(0,0)[lb]{{\SetFigFont{8}{9.6}{\rmdefault}{\mddefault}{\updefault}$pr_2$}}}
\path(600,1512)(1800,1512)
\path(600,1287)(1800,1287)
\path(600,1062)(1800,1062)
\path(600,837)(1800,837)
\path(2400,1812)(2400,612)
\path(600,12)(1800,12)
\path(1200,537)(1200,87)
\blacken\path(1185.000,147.000)(1200.000,87.000)(1215.000,147.000)(1200.000,129.000)(1185.000,147.000)
\path(1875,1212)(2325,1212)
\blacken\path(2265.000,1197.000)(2325.000,1212.000)(2265.000,1227.000)(2283.000,1212.000)(2265.000,1197.000)
\path(375,1249)(675,1512)
\path(639.771,1461.168)(675.000,1512.000)(619.994,1483.726)
\path(375,1212)(750,1287)
\path(694.107,1260.524)(750.000,1287.000)(688.223,1289.942)
\path(375,1174)(750,1062)
\path(688.217,1064.798)(750.000,1062.000)(696.802,1093.543)
\path(375,1137)(675,837)
\path(621.967,868.820)(675.000,837.000)(643.180,890.033)
\put(1875,124){\makebox(0,0)[lb]{{\SetFigFont{8}{9.6}{\rmdefault}{\mddefault}{\updefault}$A$}}}
\put(2550,762){\makebox(0,0)[lb]{{\SetFigFont{8}{9.6}{\rmdefault}{\mddefault}{\updefault}$\hat{A}$}}}
\put(150,1174){\makebox(0,0)[lb]{{\SetFigFont{8}{9.6}{\rmdefault}{\mddefault}{\updefault}$N_n$}}}
\put(0,1662){\makebox(0,0)[lb]{{\SetFigFont{8}{9.6}{\rmdefault}{\mddefault}{\updefault}$A \times \hat{A}$}}}
\path(600,1812)(1800,1812)(1800,612)
	(600,612)(600,1812)
\end{picture}
}
\end{center}
\caption{Factoring the universal cover
of an abelian variety over $k$}
\label{two}
\end{figure}
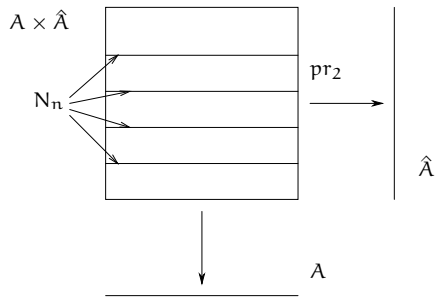

If $k$ is not separably closed, then we may apply the above
construction to $A \times_k k^s$, so $\tA  \rightarrow  A
 \times_k k^s \rightarrow  A$ gives a convenient factorization of the
universal cover. In the spirit of
Pardini, we have the following ``complementary'' factorization:  informally, although $\cL_{\chi}^{-1}$ may not be defined
over $k$, $\oplus_{\chi \text{ $n$-torsion}} \cL_{\chi}^{-1}$ is
defined over $k$ for each $n$.  We make this precise by noting that
any isogeny is dominated by $[n]$ (multiplication by $n$) for some
$n$, and that $[n]$ is defined over $k$.  Let $N_n := pr_2^* (\ker
[n])^{\red} \subset A \times \hat{A}$, where $pr_2$ is the projection
to $\hat{A}$ (see Figure~\ref{two}).  Note that if $n_1|n_2$ then we
have a canonical open and closed immersion $N_{n_1} \hookrightarrow
N_{n_2}$.  Let $\cP \rightarrow A \times \hat{A}$ be the Poincar\'e
bundle.  Then $\tilde{A} = (\underline{\Spec} \; \varinjlim
\cP|_{N_n}) \otimes_k k^s$.  In particular,
\begin{equation}\label{e:Tate}
\xymatrix{
& \tilde{A} \ar[dl] \ar[dr] \\
\underline{\Spec} \;  \varinjlim \cP|_{N_n} \ar[dr]  & & A \times_k k^s  \ar[dl]\\
& A}\end{equation}
is Cartesian.

This construction applies without change to proper $k$-schemes with
abelian fundamental group.  More generally, for any proper
geometrically connected $X/k$, this construction yields the maximal
abelian cover.

\subsection{Curves} Now consider universal covers of curves of genus $> 0$ over a field. (Curves are assumed to be finite type.)

\subsubsection{Failure of uniformization}\label{uniformization_false}Motivated by uniformization of Riemann surfaces,
one might hope that all complex (projective irreducible
nonsingular) curves of genus greater than $1$ have isomorphic
(algebraic) universal covers.  However, a short argument shows that two curves have the same
universal cover if and only if they have an isomorphic finite \'etale
cover, and another short argument shows that a curve can share such a cover
with only a countable number of other curves.  Less naively, one might
ask the same question over a countable field such as $\Qbar$.  One
motivation is the conjecture of Bogomolov and Tschinkel \cite{bt},
which states (in our language) that given two curves $C$, $C'$ of
genus greater than $1$ defined over $\Qbar$, there is a nonconstant
map $\tilde{C} \rightarrow \tilde{C'}$.  However, Mochizuki
\cite{mochizuki} (based on work of Margulis and Takeuchi) has shown
that a curve of genus $g>1$ over $\Qbar$ shares a universal cover with
only finitely many other such curves (of genus $g$).

\subsubsection{Cohomological dimension} One expects the universal cover to be simpler than the curve itself.   As a well-known example, the cohomological dimension of the universal cover is less than $2$, at least away from the characteristic of the base field (whereas for a proper curve, the cohomological dimension is at least $2$).
(Of course the case of a locally constant sheaf is simpler still.)

\begin{proposition}{Let $X$ be a smooth curve of genus $>0$ over a field $k$, and let $\tX \rightarrow X$ be the universal cover. For any integer $l$ such that $\chr(k) \nmid l$, the $l$-cohomological dimension of $\tX$ is less than or equal to $1$, i.e. for any $l$-torsion sheaf $\cF$ on the \'etale site of $\tX$, $H^i(\tX,\cF) = 0$ for $i>1.$ }\label{cohdim_UCcurve<2} \end{proposition}

(One should not expect $\tX$ to have cohomological dimension $0$ as the cohomology of sheaves supported on subschemes can register punctures in the subscheme. For instance, it is a straight forward exercise to show that for a genus $1$ curve $X$ over $\C$, the cohomological dimension of $\tX$ is $1$.)
 
For completeness, we include a sketch of a proof due to Brian Conrad, who makes no claim to originality: since $\tX$ is isomorphic to the universal cover of $X_{k^s}$, we can assume that $k$ is separably closed. One shows that $l$-torsion sheaves on $\tX$ are a direct limit of sheaves pulled back from constructible $l$-torsion sheaves on a finite \'etale cover of $X$. One then reduces to showing that for $j: U \hookrightarrow X$ an open immersion and $\cG$ a locally constant constructible $l$-torsion sheaf on $U$, $H^i(\tX, \wp^* j_! \cG) = 0,$ where $\wp$ denotes the map $\tX \rightarrow X$. Since $\tX$ is dimension $1$, only the case $i=2$ and $X$ proper needs to be considered.  Recall that $H^2(\tX, \wp^* j_! \cG) = \varinjlim H^2(Y, j_! \cG)$ where $Y$ ranges over the finite \'etale covers of $X$, and $j_! \cG$ also denotes the restriction of $j_! \cG$ to $Y$ (see for instance \cite[III \S 1 Lemma 1.16]{milnebook} whose proof references \cite[III.3]{artin}). Applying Poincar\'e duality allows us to view the maps in the direct limit as the duals of transfer maps in group cohomology $H^0(H,\cG^*_{u_0}) \rightarrow H^0(\pi_1^{\text{\'et}}(U,u_0), \cG^*_{u_0})$, where $H$ ranges over subgroups of $\pi_1^{\text{\'et}}(U,u_0)$ containing the kernel of $\pi_1^{\text{\'et}}(U,u_0) \rightarrow \pi_1^{\text{\'et}}(X,u_0)$. One shows these transfer maps are eventually $0$ as follows: let $K$ denote the kernel of $\pi_1^{\text{\'et}}(U,u_0) \rightarrow \pi_1^{\text{\'et}}(X,u_0)$. For $H$ small enough, $(\cG^*_{u_0})^{H} = (\cG^*_{u_0})^{K}$. Restricting $H$ still further to a subgroup $H'$ produces a transfer map $(\cG^*_{u_0})^{H'} \rightarrow (\cG^*_{u_0})^{H}$ which equals multiplication by the index of $H'$ in $H$. Since $\cG^*_{u_0}$ is $l$-torsion, it therefore suffices to see that we can choose $H'$ (containing $K$) in $\pi_1^{\text{\'et}}(U,u_0)$ of arbitrary $l$-power index. This follows because there are \'etale covers of $X$ of any $l$-power degree (as such a cover can be formed by pulling back the multiplication by $l^n$ map from the Jacobian).  
 
As a well-known corollary (which is simpler to prove directly),  the cohomology of a locally constant $l$-torsion sheaf $\cF$ on $X$ can be computed with profinite group cohomology: $H^i(X, \cF) = H^i(\pi_1^{\text{\'et}}(X,x_0), \cF_{x_0})$ for all $i$. (To see this, one notes that $H^1$ of a constant sheaf on $\tX$ vanishes. By Proposition \ref{cohdim_UCcurve<2}, it follows that the pullback of $\cF$ to $\tX$ has vanishing $H^i$ for all $i>0$. One then applies the Hochschild-Serre spectral sequence \cite[III \S 2 Thm 2.20]{milnebook}. Note that this corollary only requires knowing Proposition  \ref{cohdim_UCcurve<2} for $\cF$ a constant sheaf.)

Proposition \ref{cohdim_UCcurve<2}, for $\cF$ a finite constant sheaf of any order (so no assumption that the torsion order is prime to $\chr(k)$) and for $X$ an affine curve, is in \cite[Prop. 1]{serreRE} for instance. The above corollary for affine curves is in \cite[Prop. 1]{serreRE} as well. For $X$ a proper curve or affine scheme (of any dimension) and $\cF = \Z / p$ for $p = \chr(k)$, Proposition \ref{cohdim_UCcurve<2} and the Corollary are in \cite[Prop. 1.6]{Gille}. Both these references also give related results on the cohomological dimension (including $p$ cohomological dimension!) of the fundamental group of $X$, for $X$ a curve or for $X$ affine of any dimension. Also see \cite[VI \S 1]{milnebook} for related dimensional vanishing results.

\subsubsection{Picard groups} The universal covers
of elliptic curves and hyperbolic projective curves over $\C$ have
very large Picard groups, isomorphic to $(\R / \Q)^{\oplus 2}$ and
countably many copies of $\R/\Q$ respectively.
\label{relatedwarning}

\subsection{Algebraic Teichm\"uller space}
 If $g \geq 2$, then $\cm_g[n]$,
the moduli of curves with level $n$ structure, is a scheme for $n \geq 3$, and $\cm_g[n] \rightarrow \cm_g$ is finite \'etale
(where $\cm_g$ is the moduli stack of curves).  Hence $\ct :=
\tilde{\cm_g}$ is a scheme, which could be called {\em
  algebraic Teichm\"uller space}.  The {\em algebraic mapping class
  group scheme} $\upi(\cm_g)$ acts on it.  

One might hope to apply some of the methods of Teichm\"uller theory to
algebraic Teichm\"uller space.  Many ideas relating to ``profinite
Teichm\"uller theory'' appear in \cite{boggi}.  On a more analytic
note, many features of traditional Teichm\"uller theory carry over,
and have been used by dynamicists and analysts, see for example
\cite{mcmullen}.  The ``analytification'' of algebraic Teichm\"uller
space is a solenoid, and was studied for example by Markovic and \v
Sari\'c in \cite{ms}.  McMullen pointed out to us that it also yields
an interpretation of Ehrenpreis and Mazur's conjecture, that given any
two compact hyperbolic Riemann surfaces, there are finite covers of
the two surfaces that are arbitrarily close, where the meaning of
``arbitrarily close'' is not clear \cite[p.~390]{ehrenpreis}.  (Kahn
and Markovic have recently proved this conjecture using the
Weil-Petersson metric, suitably normalized, \cite{km}.)  More
precisely: a Galois type of covering of a genus $h$ curve, where the
cover has genus $g$, gives a natural
correspondence $$\xymatrix{\mathcal{X} \ar[r] \ar[d] & \cm_{g} \\
  \cm_h}$$ where the vertical map is finite \'etale.  One
might hope that the metric can be chosen on $\mathcal{M}_g$ for all
$g$ so that the pullbacks of the metrics from $\mathcal{M}_g$ and
$\mathcal{M}_h$ are the same; this would induce a pseudometric on
algebraic Teichm\"uller space. In practice, one just needs the metric to
be chosen on $\cm_g$ so that the correspondence induces a system of
metrics on $\widetilde{\cm_h}$ that converges; hence the normalization
chosen in \cite{km}.  The Ehrenpreis-Mazur conjecture asserts that
given any two points on $\mathcal{M}_h$, there exist lifts of both to
algebraic Teichmuller space whose distance is zero.

\section{The algebraic fundamental group family}
\label{s:FG}

We now construct the fundamental group family $\upi(X)$
and describe its properties.  More generally, suppose $f: Y
\rightarrow X$ is a Galois profinite-\'etale covering space with $Y$ connected.  We will
define the adjoint bundle $\Ad f: \Ad Y \rightarrow X$ to $f$, which
is a group scheme over $X$ classifying profinite-\'etale covering
spaces of $X$ whose pullback to $Y$ is Yoneda trivial. We define $\upi(X)$ as $\Ad(\tX \rightarrow X)$.

$\Ad Y$ is the quotient scheme $(Y \times_X Y) / \Aut(Y/X)$, where
$\Aut(Y/X)$ acts diagonally. The quotient is constructed by descending
$Y \times_X Y \rightarrow Y$ to an $X$-scheme, using the fact that
profinite-\'etale covering spaces are fpqc. This construction is as follows:

By Proposition \ref{trivial_Aut_bundle_Galois_covering_fiber_square}, we have the fiber square \eqref{Galois_fpqc_cover}. A descent datum on a $Y$-scheme $Z$ is equivalent to an action of $\underline{\Aut(Y/X)}_X$ on $Z$ compatible with $\mu$ in the sense that the diagram
\begin{equation}
\label{Galois_compatible_group_action_ii}
\xymatrix{\underline{\Aut(Y/X)}_X \times_X Z \ar[rr] \ar[d] && Z \ar[d] \\
\underline{\Aut(Y/X)}_X \times_X Y \ar[rr]^-{\mu} && Y}
\end{equation}
\noindent commutes. (This is the analogue of the equivalence between descent data for finite \'etale $G$ Galois covering spaces and actions of the trivial group scheme associated to $G$. The proof is identical; one notes that the diagram (\ref{Galois_compatible_group_action_ii}) is a fiber square and then proceeds in a straightforward manner. See \cite[p.~140]{neronmodels}.) 
We emphasize that for $Z$ affine over $Y$, a descent datum is easily seen to be automatically effective --- this has been a source of some confusion --- (see for instance \cite[p.~134, Thm.~4]{neronmodels}, as well as the following discussion \cite[p.~135]{neronmodels}). It follows that $\Ad Y$ exits.

\begin{definition}\label{adjoint_bundle_def} Let $f: Y
\rightarrow X$ be a Galois profinite-\'etale covering space with $Y$ connected. The {\em adjoint bundle} to $f$ is the $X$-scheme $\Ad f : \Ad Y \rightarrow X$ determined by the affine $Y$ scheme $Y \times_X Y \rightarrow Y$ and the action $\mu \times \mu.$\end{definition}

\begin{definition}\label{upi_def} Let $X$ be a scheme admitting a universal cover $\tX$. (For instance $X$ could be any connected qcqs scheme.) The {\em fundamental group family} of $X$ is defined to be $\Ad(\tX \rightarrow X),$ and is denoted $\upi(X) \rightarrow X$.  \end{definition}

$\Ad Y$ is a group scheme over $X$. The multiplication map is defined as follows: let $\Delta: Y \rightarrow Y \times_X Y$ be the diagonal map. By the same method used to construct $\Ad Y$, we can construct the $X$-scheme $(Y \times Y \times Y)/ \Aut(Y/X)$, where $\Aut(Y/X)$ acts diagonally. The map $\textrm{id} \times \Delta \times \textrm{id}: Y \times Y \times Y \rightarrow Y \times Y \times Y \times Y$ descends to an isomorphism of $X$-schemes 
\begin{equation}
\label{product_Ad_Y}
 (Y \times Y \times Y) / \Aut(Y/X) \rightarrow \Ad(Y) \times_X \Ad(Y).
 \end{equation}
The projection of $Y \times Y \times Y$ onto its first and third factors descends to a map 
\begin{equation}
\label{threeYs_to_two}
(Y \times Y \times Y)/ \Aut(Y/X) \rightarrow \Ad(Y).
\end{equation}
The multiplication map is then the inverse of isomorphism \eqref{product_Ad_Y} composed with map \eqref{threeYs_to_two}.

Heuristically, this composition law has the following description:
the
geometric points of $\Ad Y$ are equivalence classes of ordered pairs of geometric points of
$Y$ in the same fiber. Since $\Aut(Y/X)$ acts simply transitively on
the points of any fiber, such an ordered pair is equivalent to an
$\Aut(Y/X)$-invariant permutation of the corresponding fiber of $Y$
over $X$. The group law on $\Ad Y$ comes from composition of
permutations.
 
 The identity map $X \rightarrow \Ad(Y)$ is the $X$-map descended from
 the $Y$-map $\Delta$. The inverse map is induced by the map $Y
 \times_X Y \rightarrow Y \times_X Y$ which switches the two factors
 of $Y$. It is straightforward to see that these maps give $\Ad Y$ the
 structure of a group scheme.
 
 The construction of $\Ad(Y)$ implies the following:
 
\begin{proposition}{Let $Y$ be a connected profinite-\'etale Galois covering space of $X$. We have a canonical isomorphism of $Y$-schemes $\Ad(Y) \times_X Y \cong Y \times_X Y.$ Projection $Y \times_X Y \rightarrow Y$ onto the second factor of $Y$ gives an action 
  \begin{equation}
 \label{AdY_action}
 \Ad(Y) \times_X Y \rightarrow Y.
 \end{equation}
  }\end{proposition}

\begin{proposition}{Suppose $Y_1, Y_2$ are connected profinite-\'etale Galois covering spaces of $X$. An $X$-map $Y_1 \rightarrow Y_2$  gives rise to a morphism of group schemes $\Ad(Y_1) \rightarrow \Ad(Y_2)$. Furthermore, the map $\Ad(Y_1) \rightarrow \Ad(Y_2)$ is independent of the choice of $Y_1 \rightarrow Y_2$. }\label{Ad_functoriality_fixed_base}\end{proposition}
 
\begin{proof}
 Choose a map $g: Y_1 \rightarrow Y_2$ over $X$. By Proposition \ref{trivial_Aut_bundle_Galois_covering_fiber_square}, we have an isomorphism $Y_2 \times_X Y_2 \rightarrow \underline{\Aut(Y_2/X)}_{Y_2}$ defined over $Y_2$.  Pulling this isomorphism back by $g$ gives an isomorphism $Y_2 \times_X Y_1 \rightarrow \underline{\Aut(Y_2/X)}_{Y_1}$ defined over $Y_1$. The map $g \times id: Y_1 \times_X Y_1 \rightarrow Y_2 \times_X Y_1$ therefore gives rise to a $Y_1$ map $\underline{\Aut(Y_1/X)}_{Y_1} \rightarrow
 \underline{\Aut(Y_2/X)}_{Y_1}$. This map corresponds to a continuous map of
 topological spaces $\Aut(g):\Aut(Y_1/X) \rightarrow \Aut(Y_2/X)$ by
 Proposition \ref{Yoneda_triv_means_triv_bundle}.
 
 It follows from the construction of the isomorphism of Proposition \ref{trivial_Aut_bundle_Galois_covering_fiber_square} (which is really given in Proposition \ref{Yoneda_triv_means_triv_bundle}) that for any $a \in \Aut(Y_1/X)$ the diagram:
  \begin{equation}
 \label{AutY_1_to_AutY_2}
 \xymatrix{ Y_1 \ar[rr]^{a} \ar[d]^{g}&& Y_1 \ar[d]^{g} \\
 Y_2 \ar[rr]^{\Aut(g)(a)}&& Y_2}
 \end{equation}
   \noindent commutes. 
  
Since $g: Y_1 \rightarrow Y_2$ is a profinite-\'etale covering space and in particular an fpqc cover, $g^*: \Maps(Y_2, -) \rightarrow  \Maps(Y_1, -)$ is an injection. By \eqref{AutY_1_to_AutY_2}, the maps $ \Aut(g)(a_2) \circ \Aut(g) (a_1)$ and $\Aut(g)(a_2 \circ a_1)$ have the same image under $g^*$. Thus, $\Aut(g)$ is a continuous group homomorphism. 
 
 It follows that the map $g \times g: Y_1 \times Y_1 \rightarrow Y_2 \times Y_2$ determines a map $(Y_1 \times Y_1)/\Aut(Y_1/X) \rightarrow (Y_2 \times Y_2)/\Aut(Y_2/X).$ It is straightforward to see this is a map of group schemes $\Ad(Y_1) \rightarrow \Ad(Y_2).$
 
 Given two maps of $X$-schemes $g_1,g_2: Y_1 \rightarrow Y_2,$ we have a map $(g_1, g_2): Y_1 \rightarrow Y_2 \times_X Y_2.$ Since $Y_2 \times_X Y_2 \rightarrow Y_2$ is Yoneda trivial with distinguished sections $\Aut(Y_2/X),$ we have $a \in \Aut(Y_2/X)$ such that $a \circ g_1= g_2$. It follows that $g_1$ and $g_2$ determine the same map $\Ad(Y_1) \rightarrow \Ad(Y_2).$
\end{proof}
 
\begin{corollary}{If a universal cover of $X$ exists (e.g.\ if
   $X$ is connected and qcqs,
   Proposition~\ref{initial_pointed_covering_space}), $\upi(X)$ is unique up to
distinguished isomorphism, and in particular is independent of choice
of universal cover.}  \label{pi1unique}\end{corollary}
 
\begin{theorem}{There is a canonical homeomorphism between
the underlying topological group of the fiber of
$ \upi(X) \rightarrow X$ over a geometric point $x_0: \Spec \Omega \rightarrow X$ and 
the \'etale (pointed) fundamental group $\pi_1(X, x_0)$.} \label{etalefglink}\end{theorem}
 
\begin{proof} Let $Y \rightarrow X$ be a finite \'etale Galois covering space
with $Y$ connected. We have a canonical action of $X$-schemes $\upi(X)
\times_X Y \rightarrow Y$ as follows: choose a universal cover $p: \tX
\rightarrow X$ and a map $\tX \rightarrow Y$ over $X$. By Proposition
\ref{Ad_functoriality_fixed_base}, we have a canonical map $\upi(X)
\rightarrow Ad (Y).$ Composing with the canonical action $\Ad(Y)
\times_X Y \rightarrow Y$ given by \eqref{AdY_action} gives the action
$\upi(X) \times_X Y \rightarrow Y$.

Let $\fT \upi(X,x_0)$ be the topological group underlying the fiber of
$ \upi(X) \rightarrow X$ above $x_0.$ Let $\cF_{x_0}$ be the fiber functor over $x_0$. The action $\upi(X) \times_X Y \rightarrow Y$ shows that $\cF_{x_0}$ induces a functor from finite, \'etale, connected, Galois covering spaces to continuous, finite, transitive, symmetric $\fT \upi(X,x_0)$-sets. (A symmetric transitive $G$-set for a group $G$ is defined to mean a $G$-set isomorphic to the set of cosets of a normal subgroup. Equivalently, a symmetric transitive $G$-set is a set with a transitive action of $G$ such that for any two elements of the set, there is a morphism of $G$-sets taking the first to the second.)

Since $\pi_1(X, x_0)$ is characterized by the fact that $\cF_{x_0}$ induces an equivalence of categories from finite, \'etale, connected, Galois covering spaces to continuous, finite, transitive, symmetric $\pi_1(X,x_0)$-sets, it suffices to show that $\cF_{x_0}$ viewed as a functor to $\fT \upi(X,x_0)$-sets as in the previous paragraph is an equivalence of categories.
By fpqc descent, pull-back by $p$, denoted $p^*$, is an equivalence of categories
 from affine $X$-schemes to affine $\tX$-schemes with descent
data. Because $\tX$ trivializes any finite, \'etale $X$-scheme, it is
straightforward to see that $p^*$ gives an equivalence of categories
from finite, \'etale, covering spaces of $X$ to trivial, finite,
\'etale covering spaces of $\tX$ equipped with an action of
$\underline{\Aut(\tX/X)}.$ It follows from Proposition \ref{Yoneda_triv_means_triv_bundle} that taking the topological space underlying the fiber over a geometric point of $\tX$ is an equivalence of categories from trivial, finite,
\'etale covering spaces of $\tX$ equipped with an action of
$\underline{\Aut(\tX/X)}$ to continuous, finite, transitive, symmetric $\Aut(\tX/X)$-sets. Forgetting the choice of geometric point of $\tX$ shows that  $\cF_{x_0}$ is
an equivalence from the category of finite, \'etale, connected, Galois covering spaces of $X$
to continuous, finite, transitive, symmetric $\fT \upi(X,x_0)$-sets.
\end{proof}
 
 The remainder of this section is devoted to examples and properties of the fundamental group family.
 
 \subsection{Group schemes} We continue the discussion of \S
 \ref{UCgrp_sch_b_point} to obtain the algebraic version of the fact
 that if $X$ is a topological group with identity $e$, there is a
 canonical exact sequence
$$\xymatrix{ 
0
\ar[r] & \pi_1(X,e) \ar[r] & \tX \ar[r] & X \ar[r] & 0.}
$$

\begin{theorem}{If $X$ is a connected group scheme finite type over
  an algebraically closed field $k$ such that $\tX$ is a group scheme (e.g.\
  if $\chr k = 0$ or $X$ is proper, Thm.~ \ref{UCgroupscheme}), then
  the kernel of the morphism $\tX \rightarrow X$ is naturally
  isomorphic to $\upi(X,e)$ (as group schemes).} \label{groupES}\end{theorem}

\begin{proof}
Let $G$ be the (scheme-theoretic) kernel of $p: \tX \rightarrow X$. Restricting the $X$-action $$\upi(X) \times_X \tX \rightarrow \tX$$
to $e$ yields a $k$-action 
\begin{equation}
\label{january}
\upi(X,e) \times G \rightarrow G.
\end{equation}
Evaluating (\ref{january}) on $\tilde{e} \hookrightarrow G$ yields
an isomorphism $\gamma: \upi(X,e) \rightarrow G$
(using Theorem~\ref{etalefglink}). We check that $\gamma$ respects the group scheme
structures on both sides.  It suffices to check that the multiplication
maps are the same. Let $m_{\upi(X,e)}$ and $m_G$ be the multiplication maps for $\upi(X,e)$ and $G$ respectively. The diagram
$$
\xymatrix{
\upi(X,e) \times \upi(X,e) \ar[rr]^{\quad m_{\upi(X,e)}} \ar[d]^{\text{id} \times \gamma} &   & \upi(X,e) \ar[d]^{\gamma} \\
\upi(X,e) \times G  \ar[rr]^{\text{\eqref{january}}}\ar[d]^{\gamma \times \text{id}} &  & G \ar[d]    \\
G \times G \ar[rr]^{m_G} &  &  G 
}
$$
commutes. (The upper square commutes because \eqref{january} is a group action. The lower square commutes because monodromy commutes with morphisms of profinite-\'etale covering spaces. In particular, right multiplication in $\tilde{X}$ by any geometric point of $G$ commutes with the monodromy action $\upi(X) \times_X \tX \rightarrow \tX$.) This gives the result.
\end{proof}

\subsection{Examples} \label{pi1examples}
We now describe the fundamental group family in a number of
cases.

\subsection{The absolute Galois group scheme}\label{fund_grp_fmly_G_Q}
We give four descriptions of the {\em absolute Galois group scheme} $\underline{\Gal}(\Q):=
\upi(\Spec \Q)$, or equivalently, we describe the corresponding Hopf algebra. As  $\underline{\Gal}(\Q)$ does not depend on the choice of the algebraic closure $\Qbar$ (Prop.~\ref{pi1unique}), we do not call it $\underline{\Gal}(\Qbar/\Q)$. {\em Notational Caution:}  $\underline{\Gal}(\Q)$ is not the trivial group scheme corresponding
to $\Gal(\Qbar/\Q)$, which would be denoted $\underline{\Gal(\Qbar/\Q)}$ (Example \ref{trivialprofinitegroup}).

{\em 1) By definition.} The Hopf algebra
consists of those elements of $\Qbar \otimes_{\Q} \Qbar$ that are
invariant under the diagonal action of the Galois group
$\Gal(\Qbar/\Q)$. The coidentity map sends $a \otimes b$ to $ab$. The coinverse map is given by the involution $a
\otimes b \mapsto b \otimes a$.  The comultiplication map has the following description: $id \otimes \De \otimes id$ gives a map
  $\otimes^4 \Qbar \rightarrow  \otimes^3 \Qbar$ which descends to an isomorphism
  $\otimes^2 ((\Qbar \otimes \Qbar)^{\Gal(\Qbar/ \Q)}) \rightarrow (\otimes^3 \Qbar)^{\Gal (\Qbar/ \Q)}$, where all
  actions of $\Gal(\Qbar/ \Q)$ are diagonal. The comultiplication map can therefore be
  viewed as a map $(\otimes^2 \Qbar)^{\Gal(\Qbar/\Q)} \rightarrow (\otimes^3 \Qbar)^{\Gal (\Qbar/ \Q)}$ and this map
  is the inclusion onto the first and third factors.

{\em 2) As an algebra of continuous maps.}  The Hopf algebra consists of continuous
maps $f: \Gal(\Qbar / \Q) \rightarrow \Qbar$ such that
\begin{equation} \label{compatibility}
\xymatrix{
\Gal(\Qbar / \Q) \ar[d]^{\si} \ar[r]^{\quad f} & \Qbar \ar[d]^{\si} \\
\Gal(\Qbar / \Q)  \ar[r]^{\quad f} & \Qbar \\
}
\end{equation}
commutes for all $\si \in \Gal(\Qbar/\Q)$, 
where the left vertical arrow is conjugation, and the
right vertical arrow is the Galois action.
Note that these maps form an algebra.
The coinverse of $f$ is given by the composition
$\xymatrix{ \Gal(\Qbar / \Q)  \ar[r]^i & \Gal(\Qbar / \Q) \ar[r]^{\quad f} & \Qbar}$, where $i$ is the inverse in $\Gal(\Qbar / \Q)$. Comultiplication applied to $f$ is given by the composition
$$\xymatrix{ \Gal(\Qbar / \Q)   \times \Gal(\Qbar / \Q)  
  \ar[r]^{\quad \quad m} & \Gal(\Qbar / \Q) \ar[r]^{\quad f} &
  \Qbar},$$ using the isomorphism 
$$\Maps_{\rm{cts}}(\Gal(\Qbar/\Q)
\times \Gal(\Qbar/\Q), \Qbar) \cong \Maps_{\rm{cts}}(\Gal(\Qbar/\Q),
\Qbar ) \otimes \Maps_{\rm{cts}}(\Gal(\Qbar/\Q, \Qbar).$$ (A similar argument was used to construct the trivial profinite group scheme
in Example \ref{trivialprofinitegroup}. The similarity comes from the isomorphism of $\upi \times_X \tX$ with $\underline{\Aut(\tX/X)}_{\tX}$.)

{\em 3) Via finite-dimensional representations.}  By interpreting
\eqref{compatibility} as ``twisted class functions,'' we can describe
the absolute Galois Hopf algebra in terms of the irreducible
continuous representations of $\Gal(\Qbar/\Q)$ over $\Q$.  More precisely, we give a
basis of the Hopf algebra where comultiplication and coinversion are
block-diagonal, and this basis is described in terms of
$\Q$-representations of $\Gal(\Qbar/\Q)$.

Given a finite group $G$ and a representation $V$ of $G$ over a field
$k$, the natural map $G \rightarrow V\otimes V^*$ induces a map  \[(V\otimes V^*)^*\rightarrow \Maps(G,
k),\] where $V^*$ denotes the dual vector space. For simplicity, assume that $k$ is a subfield of $\C.$ When $k$ is algebraically closed, Schur orthogonality gives
that \[ \Maps(G,k) \cong \oplus_{V \in I} (V \otimes V^*)^*,\] where
$I$ is the set of isomorphism classes of irreducible representations
of $G$. It follows that \[ \Maps_{\rm{cts}}(\Gal(\Qbar/\Q), \Qbar)
\cong \oplus_{G \in Q} \oplus_{V \in I_G} (V \otimes V^*)^*\] where
$Q$ is the set of finite quotients of $\Gal(\Qbar/\Q)$, and for any
$G$ in $Q$, $I_G$ is the set of isomorphism classes of irreducible,
faithful representations of $G$ over $\Qbar.$

 $\Gal(\Qbar/ \Q )$ acts on $\Maps_{\rm{cts}}(\Gal(\Qbar/\Q), \Qbar)$ via $(\sigma f) (\sigma') = \sigma (f(\sigma^{-1} \sigma' \sigma ))$, where $f: \Gal(\Qbar/\Q) \rightarrow \Qbar$ is a continuous function and $\sigma$, $\sigma'$ are in $\Gal(\Qbar/\Q).$ The set of fixed points is the Hopf algebra we wish to describe. The elements of this Hopf algebra could reasonably be called ``twisted class functions''. Note that we have a $\Q$-linear projection from $\Maps_{\rm{cts}}(\Gal(\Qbar/ \Q), \Qbar)$ to our Hopf algebra given by averaging the finite orbit of a function. 
 
 Let $G$ be a finite quotient of $\Gal(\Qbar/\Q).$ $\Gal(\Qbar/\Q)$
 acts on the irreducible, faithful $\Qbar$-representations of $G$ by
 tensor product, namely, $\sigma(V) = \Qbar \otimes_{\Qbar} V,$ where
 the map $\Qbar \rightarrow \Qbar$ in the tensor product is $\sigma.$
 The orbits of $I_G$ under this action are in bijection with the
 irreducible, faithful $\Q$-representations of $G$. This bijection sends an irreducible, faithful
 $\Qbar$-representation $V$ to the isomorphism class of
 $\Q$-representation $W_V$ such that \[ \oplus_{W \in O_V} W \cong W_V
 \otimes \Qbar\] where $O_V$ is the (finite) orbit of $V$ under the
 action of $\Gal(\Qbar/\Q)$.
 
 For any irreducible, faithful $\Qbar$-representation $V$ of $G$, $\oplus_{W \in O_V} (W \otimes W^*)^*$ is an invariant subspace of $\Maps_{\rm{cts}}(\Gal(\Qbar/\Q), \Qbar)$ under the action of $\Gal(\Qbar/\Q)$.  It follows that our Hopf algebra is isomorphic to \[ \oplus_{G \in Q} \oplus_{V \in \overline{I}_G } (\oplus_{W \in O_V}(W \otimes W^*)^*)^{\Gal(\Qbar/\Q)} \] where  $\overline{I}_G$ is the set of orbits of $I_G$ under $\Gal(\Qbar/\Q).$
 
 The natural map $\oplus_{W \in O_V}(W \otimes W^*)^*\rightarrow
 \Maps_{\rm{cts}}(\Gal(\Qbar/\Q), \Qbar)$ factors through the natural
 map $(W_V \otimes W_V^* \otimes \Qbar)^*\rightarrow
 \Maps_{\rm{cts}}(\Gal(\Qbar/\Q), \Qbar)$. Note that there is a
 compatible $\Gal(\Qbar/\Q)$-action on $(W_V \otimes W_V^* \otimes
 \Qbar)^*$. Note that the map $(W_V \otimes W_V^* \otimes
 \Qbar)^*\rightarrow \Maps_{\rm{cts}}(\Gal(\Qbar/\Q), \Qbar)$ is not
 injective.  Let the image of $(W_V \otimes W_V^* \otimes \Qbar)^*$ in
 $\Maps_{\rm{cts}}(\Qbar/\Q), \Qbar)$ be $\mathcal{F}(W_V)$.
 
 Let $\overline{I}_{{\Gal(\Qbar/\Q)}}$ be the set of isomorphism classes of continuous irreducible $\Q$-representations of $\Gal(\Qbar/\Q)$. Our Hopf algebra is isomorphic to \[ \oplus_{\overline{I}_{{\Gal(\Qbar/\Q)}}} \mathcal{F}(W_V)^{\Gal(\Qbar/\Q)}.\] The subspaces $\mathcal{F}(W_V)^{\Gal(\Qbar/\Q)}$ are invariant under comultiplication and coinversion because comultiplication and coinversion are induced from comultiplication and coinversion on $GL(W_V \otimes \Qbar)$. The multiplication is not diagonal; it comes from tensor products of representations and therefore involves the decomposition into irreducible representations of the tensor product of two irreducible representations.
  
{\em 4) Points of the absolute Galois group scheme.}  Let $K
\rightarrow L$ be a finite Galois extension of fields with Galois
group $G$.  The points and group scheme structure of the adjoint bundle $\Ad(L/K) := \Ad(\Spec L \rightarrow \Spec K)$ can be identified as follows: as in part 2) of this example, the
ring of functions of $\Ad(L/K)$ is the ring of functions $f: G \rightarrow
L$ such that for all $g,h$ in $G$, $f(hgh^{-1}) = h f(g)$. Thus, the
points of $\Ad(L/K)$ are in bijection with conjugacy classes of
$G$. Specifically, let $S$ be a set of representatives of the
conjugacy classes of $G$. For any element $g$ of $G$, let $C_g$ be the
centralizer of $g$. Then
$\Ad(L/K) = \coprod_{c \in S} \Spec L^{C_{c}}$. 

The group law on $\Ad(L/K)$ therefore corresponds to a map $\coprod_{a,b\in S}
\Spec (L^{C_{a}} \otimes L^{C_{b}}) \rightarrow \coprod_{c \in S}
\Spec L^{C_{c}}$.  Note that $\Spec (L^{C_{a}} \otimes L^{C_{b}}) =
\coprod_{g \in S_{a,b}} \Spec (L^{C_a} (g L^{C_b}))$, where $S_{a,b}$
is a set of double coset representatives for $(C_a, C_b)$ in $G$,
i.e.\ $G = \coprod_{g\in S_{a,b}} C_a g C_b$, and $L^{C_a} (g
L^{C_b})$ is the subfield of $L$ generated by $L^{C_a}$ and $g
L^{C_b}$.  (In particular, the points of $\Spec (L^{C_{a}} \otimes
L^{C_{b}})$ are in bijective correspondence with $S_{a,b}$.) Noting
that $L^{C_a} (g L^{C_b})= L^{C_a \cap gC_b g^{-1}}$, we have that
the comultiplication on $\Ad$ is a map
\begin{equation}
\label{comult_Ad(Galois_ext)}
 \prod_{c \in S} L^{C_c} \rightarrow \prod_{a,b \in S} \prod_{g \in S_{a,b}} L^{C_a \cap C_{gbg^{-1}}}
 \end{equation} 

Comultiplication is described as follows: $L^{C_c} \rightarrow L^{C_a \cap C_{gbg^{-1}}}$is the $0$ map if $c$ is not contained in the set $R_{a,b} = \{ g_1 a g_1^{-1} g_2 b g_2^{-1} | g_1, g_2 \in G \}$. Otherwise, there exists $g'$ in $G$ such that $g' c {g'}^{-1} = a g b g^{-1}$. The map $L^{C_c} \rightarrow L^{C_a \cap C_{gbg^{-1}}}$ is then the composite \[ L^{C_c} \stackrel{g'}{\rightarrow} L^{C_{g' c {g'}^{-1}}}=L^{C_{a g b g^{-1}}} \hookrightarrow L^{C_a \cap C_{gbg^{-1}}}.\] 

Note that $R_{a,b}$ is a union of conjugacy classes, and these conjugacy classes are in bijection with $S_{a,b}$, just like the points of $\Spec (L^{C_{a}} \otimes L^{C_{b}})$.

This description is explicit; the reader could easily write down the comultiplication map for the $S_3$ Galois extension $\Q \rightarrow \Q(2^{1/3}, \omega)$, where $\omega$ is a primitive third root of unity.

We obtain the following description of $\underline{Gal}(\Q) = \upi(\Spec \Q)$: replace the
products in \eqref{comult_Ad(Galois_ext)} by the subset of the
products consisting of continuous functions. The map (\ref{comult_Ad(Galois_ext)}) restricts to the comultiplication map between these function spaces.

\subsubsection{Residue fields of  $\underline{Gal}(\Q)$}
\label{pointsofGal} Note that the points of $\upi(\Spec \Q)$
correspond to conjugacy classes in $\Gal(\Qbar/\Q)$, and their residue
fields are the fixed fields of the centralizers. Although any two commuting elements of $\Gal(\Qbar/\Q)$ are contained in a copy of $\hat{\Z}$ or $\Z/2$ (\cite{Geyer}), centralizers are not necessarily even abelian. Indeed a ``folklore" theorem told to us by Florian Pop says that every countably generated group of cohomological dimension $1$ is a subgroup of $\Gal(\Qbar/\Q).$ In particular, for two distinct primes $l_1$ and $l_2$ and any action of $\hat{\Z}_{l_1}$ on $\hat{\Z}_{l_2}$, $\hat{\Z}_{l_2} \rtimes \hat{\Z}_{l_1}$ is a subgroup of $\Gal(\Qbar/\Q).$ If we choose $l_1$ and $l_2$ such that the $l_1$ Sylow subgroup of $\hat{\Z}_{l_2}^*$ is non-trivial, we may choose a non-trivial action of $\hat{\Z}_{l_1}$ on $\hat{\Z}_{l_2}$, yielding a non-abelian group $\hat{\Z}_{l_2} \rtimes \hat{\Z}_{l_1}$. The center of $\hat{\Z}_{l_2} \rtimes \hat{\Z}_{l_1}$ is non-trivial. It follows that the non-abelian group $\hat{\Z}_{l_2} \rtimes \hat{\Z}_{l_1}$ is contained in a centralizer of  $\Gal(\Qbar/\Q)$.

\subsection{Finite fields $\F_q$}
Parts 1), 2) and 4) of Example \ref{fund_grp_fmly_G_Q} apply to any field $k$, where $\Qbar$ is
replaced by $k^s$.  In the case of a finite field, the Galois group is
abelian, so the compatibility condition \eqref{compatibility}
translates to the requirement that a continuous map  $\Gal(\overline{\F_q} / \F_q)
\rightarrow \overline{\F_q}$ have image contained in $\F_q.$
  Hence, $\upi(\Spec \F_q)$ is the trivial
profinite group scheme $\underline{\hat{\Z}}$ over $\F_q$ (see Example \ref{trivialprofinitegroup}).

\subsection{$\G_m$ over an algebraically closed field $k$ of characteristic $0$}
Note that $\Gamma(\tilde{\G}_m \times_{\G_m} \tilde{\G}_m)$ can be
interpreted as the ring $k[u_1^{\Q}, u_2^{\Q}]$ subject to $u_1^n =
u_2^n$ for $n$ in $ \Z$ (but not for general $n \in \Q$).  Thus
$\upi(\G_m) = (k[u_1^{\Q}] \otimes_{k[t^\Z]} k[u_2^{\Q}])^{\Aut (
  k[t^{\Q}]/ k[t^{\Z}])}$. 
The automorphisms of $k[t^{\Q}]/ k[t^{\Z}]$ involve sending $t^{1/n}$ to
$\zeta_n t^{1/n}$, where $\zeta_n$ is an $n$th root of unity, and all
the $\zeta_n$ are chosen compatibly.  Hence the invariants may be
identified with $k[t^\Z][\mu_{\infty}]$ where $t_{n!} =
(u_1/u_2)^{1/n!}$. Thus we recognize the fundamental group scheme as  $\underline{\hat{\Z}}$
(Example \ref{Zhat}). The action of $\upi(\G_m)$ on $ \tilde{\G}_m$ is given by
$$
\xymatrix{
k[t^{\Q}] \ar[r] & k[t^{\Q}] \otimes_{k[t,t^{-1}]} k[t^{\Z}, t_1, \dots]/
(t_1-1, t_{n!}^n-t_{(n-1)!})}
$$
with $t^{1/n!} \mapsto t_{n!} t^{1/n!}.$
Notice that we get a natural exact sequence of group schemes over $k$
$$
\xymatrix{
0 \ar[r] & \underline{\hat{\Z}} \ar[r] & \tilde{\G}_m \ar[r] & \G_m \ar[r] & 0,}
$$
which is Theorem~\ref{groupES} in this setting.

In analogy with Galois theory, we have:

\begin{proposition}
{\label{XYZ}Suppose $f: X \rightarrow Y$, $g: Y \rightarrow Z$,
and $h = g \circ f$ are profinite-\'etale
covering spaces with $X$, $Y$, and $Z$ connected.
\begin{enumerate}
\item[(a)] If $h$ is Galois, then $f$ is Galois. 
There is a natural closed immersion of group schemes on $Y$
$\Ad(X/Y) \hookrightarrow g^* \Ad(X/Z)$.
\item[(b)]  If furthermore $g$ is Galois, then we have a
natural surjection
$\Ad(X/Z) \rightarrow \Ad(Y/Z)$ of group schemes over $Z$.
The kernel, which we denote $\Ad_Z(X/Y)$, is a group scheme over
$Z$
$$
1 \rightarrow \Ad_Z(X/Y) \rightarrow \Ad(X/Z) \rightarrow \Ad(Y/Z) \rightarrow 1
$$
and upon pulling this sequence back by $g$, we obtain an isomorphism
$g^* \Ad_Z(X/Y) \cong \Ad(X/Y)$ commuting with the inclusion of (a):
$$
\xymatrix{
g^* \Ad_Z(X/Y) \ar@{<->}[d]^{\sim} \ar@{^(->}[r] & g^*\Ad(X/Z)  \\
\Ad(X/Y) \ar@{^(->}[ur]
}
$$
\item[(c)] If furthermore $\Aut(X/Y)$ is abelian, then we have an action of $\Ad(Y/Z)$ on $\Ad_Z(X/Y)$, which when pulled back to $X$ is the action $$\underline{\Aut(Y/Z)}_X \times_X \underline{\Aut(X/Y)}_X \rightarrow  \underline{\Aut(X/Y)}_X$$ arising from the short exact sequence with abelian kernel $$1 \rightarrow \Aut(X/Y) \rightarrow \Aut(X/Z) \rightarrow \Aut(Y/Z) \rightarrow 1.$$ (Recall that to any short exact sequence of groups $1 \rightarrow A \rightarrow B \rightarrow C \rightarrow 1$ with $A$ abelian, $C$ acts on $A$ by $c (a) := bab^{-1}$ where $b$ is any element of $B$ mapping to $c$.)
\end{enumerate}
}\end{proposition}

We omit the proof, which is a straightforward verification.

\subsection{$\G_m$ over a field $k$ of characteristic $0$}
\label{Gmgeneral}We now extend the previous example to 
an arbitrary field of characteristic $0$.  The universal cover
of $\Spec k[t^{\Z}]$ is $\Spec \kbar[t^{\Q}]$.   

 Consider the diagram
$$
\xymatrix{
& \Spec \kbar[t^{\Q}] \ar[dl]_{j^*d^*\uGal ( k)}^i \ar[dr]^f \\
\Spec k[t^{\Q}]   \ar[dr]^j_{\text{not Galois}} & \square &  \Spec \kbar[t^{\Z}] \ar[dl]_{d^*\uGal (k)}^g \ar[dr]^{\text{not profinite-\'etale}} \\
   & \Spec k[t^{\Z}]  \ar[dr]^d_{\text{not profinite-\'etale}} & \square & \Spec \kbar \ar[dl]^{\uGal ( k)} \\
  & & \Spec k
}
$$
in which both squares are Cartesian.
All but the two indicated morphisms are profinite-\'etale. By base change from $\Spec \kbar \rightarrow \Spec k$, we see that
each of the top-right-to-bottom-left morphisms is Galois with
adjoint bundle given by the pullback of $\uGal (k)$. 
(Note:
$\Spec k[t^{\Q}] \rightarrow \Spec k[t^{\Z}]$ is {\em not}
Galois in general.)
By Proposition~\ref{XYZ}(b), with $f$ and $g$ used in the same sense,
we have an exact sequence of group schemes on $\G_m = \Spec k[t^{\Z}]$:
\begin{equation}
\label{Tate_module_exact_sequence_G_m}
\xymatrix{ 1 \ar[r] & T \ar[r] & \upi(\G_m) \ar[r] & d^* \uGal (k)
\ar[r] & 1.}
\end{equation}

Since $T$ is abelian, we have an action of $d^* \uGal( k)$ on $T$ by Proposition~\ref{XYZ}(c). 

By Proposition~\ref{XYZ}(a) applied to $\Spec \kbar[t^{\Q}] \rightarrow \Spec k[t^{\Q}]  \rightarrow  \Spec k[t^{\Z}]$, the exact sequence (\ref{Tate_module_exact_sequence_G_m}) is split when pulled back to $\Spec k[t^\Q].$

(\ref{Tate_module_exact_sequence_G_m}) is independent of the choice of algebraic closure by Corollary \ref{pi1unique}.
If we examine this exact sequence over the geometric point 
$\tilde{e} = \Spec \kbar$ mapping to the identity in $\G_m$, we obtain
\begin{equation}\label{geom_arith_pi_1_G_m}
\xymatrix{ 1 \ar[r] & \underline{\hat{\Z}} \ar[r] & \upi(\G_m,
  \tilde{e}) \ar[r] & \underline{\Gal(\kbar/ k)} \ar[r] & 1}
\end{equation}
inducing a group scheme action 
\begin{equation}\label{cyclotomic}
\underline{\Gal(\kbar/k)} \times \underline{\hat{\Z}} 
\rightarrow
\underline{\hat{\Z}}.
\end{equation}
If $k=\Q$, the underlying topological space of (\ref{geom_arith_pi_1_G_m}) (forgetting the scheme
structure) is the classical exact sequence (e.g.\ \cite[p.~77]{oort})
$$
\xymatrix{ 0 \ar[r] & \hat{\Z} \ar[r] & \pi^{et}_1(\proj^1_{\Q} - \{ 0, \infty \},1)
\ar[r] & \Gal(\Qbar/ \Q) \ar[r] & 0}$$
and the representation \eqref{cyclotomic} is a schematic version of the 
cyclotomic representation  $\rho:  \Gal(\Qbar/\Q) \rightarrow \Aut(\hat{\Z})$.   
  
\subsection{Abelian varieties}\label{av} The analogous argument holds for
an abelian variety $A$ over any field $k$. Using the diagram
$$
\xymatrix{
& \tilde{A}\ar[dl]_{j^*d^*\uGal(k)} \ar[dr] \\
\underline{\Spec} \; \varinjlim \cP|_{N_n}   \ar[dr]^j_{\text{not Galois}} & \square &  A \times_k k^s\ar[dl]_{d^*\uGal(k)} \ar[dr]^{\text{not profinite-\'etale}} \\
   & A   \ar[dr]^d_{\text{not profinite-\'etale}} & \square & \Spec k^s \ar[dl]^{\underline{\Gal}(k)} \\
  & & \Spec k
}
$$
we obtain an exact sequence of group schemes over $A$
$$
1 \rightarrow T \rightarrow \upi(A) \rightarrow d^*(\uGal(k)) \rightarrow 1
$$ 
inducing a canonical group scheme action
\begin{equation}
\label{action}
d^* \uGal(k) \times T \rightarrow T.
\end{equation}
Upon base change to the geometric point $\tilde{e} = \Spec k^s$, we obtain
$$
1 \rightarrow \underline{T'} \rightarrow \upi(A, \tilde{e}) \rightarrow
\underline{\Gal(k^s/k)} \rightarrow 0$$

(where $T' \cong \hat{\Z}^{2g}$ if $\chr k=0$, and the obvious variation in positive
characteristic), and the group action \eqref{action}
becomes the classical Galois action on the Tate module.  

More generally, for any geometrically connected $k$-variety
with a $k$-point $p$, the same argument gives a schematic version of \cite[
Exp.\ X.2, Cor.~2.2]{sga1}.

\subsection{Algebraic $K(\pi,1)$'s and elliptic curves}\label{K(pi,1)_elliptic_curve} 
(In this discussion, note that the phrase $K(\pi,1)$ has a well-established meaning in arithmetic
geometry.  We are discussing {\em different} possible analogies of this topological notion, and hope
no confusion will result.)
We suggest (naively) a direction in which to search for alternate definitions of ``trivial covering space'' and ``covering space'' to use in the procedure to produce a fundamental group family described in the introduction. For simplicity, we restrict our attention to schemes over a given number field $k$. Homomorphisms between \'etale fundamental groups are also assumed to respect the structure map to $\Gal(\kbar/k)$ up to inner automorphism. (The condition ``up to inner automorphism'' comes from ambiguity of the choice of base point, which is not important for this example, but see \cite{SzLuminy} for a careful treatment.)

The question ``what is a loop up to homotopy?'' naturally leads to the question ``which spaces are determined by their loops up to homotopy?'' When a ``loop up to homotopy'' is considered to be an element of the \'etale fundamental group, a well-known answer to the latter question was conjectured by Grothendieck: in \cite{Groth_Brief_Faltings}, Grothendieck conjectures the existence of a subcategory of ``anabelian'' schemes, including hyperbolic curves over $k$, $\Spec k$, moduli spaces of curves, and total spaces of fibrations with base and fiber anabelian, which are determined by their \'etale fundamental groups. These conjectures can be viewed as follows: algebraic maps are so rigid that homotopies do not deform one into another. From this point of view, a $K(\pi,1)$ in algebraic geometry could be defined as a variety $X$ such that $\operatorname{Mor}(Y,X) = \operatorname{Hom}(\pi_1(Y),\pi_1(X)),$ for all reasonable connected schemes $Y$. (Again, more care should be taken with base points, but this is not important here.)  For this example, define a  scheme to be a {\em $K(\pi,1)$ with respect to the \'etale fundamental group} in this manner, where $\pi_1$ is taken to be the \'etale fundamental group . In other words, ``anabelian schemes'' are algebraic geometry's $K(\pi,1)$'s with respect to the \'etale fundamental group. (Some references on the anabelian conjectures are \cite{Groth_Brief_Faltings,
 Esquisse, coh_num_fields, pop, SzLuminy}. For context, note that one could define a scheme $X$ to be a $K(\pi,1)$ if $X$ has the \'etale homotopy type of $B\pi_1^{\text{\'et}}(X)$ as in \cite{am}, but that this is not the definition we are using.) 
 
 From the above list, we see that Grothendieck conjectures that many familiar $K(\pi,1)$'s from topology are also $K(\pi,1)$'s with respect to the \'etale fundamental group in algebraic geometry, but that elliptic curves and abelian varieties are notably omitted from this list. There are many straightforward reasons for the necessity of these omissions, but since we are interested in what a loop up to homotopy should be, we consider the following two, both pointed out to us by Jordan Ellenberg.
 
By a theorem of Faltings,  $\operatorname{Mor}(Y,X) \otimes \hat{\Z} \rightarrow \operatorname{Hom}(\pi_1^{\text{\'et}}(Y),\pi_1^{\text{\'et}}(X))$ is an isomorphism for two elliptic curves $X$ and $Y$. In particular, although $\operatorname{Mor}(Y,X) \rightarrow \operatorname{Hom}(\pi_1^{\text{\'et}}(Y),\pi_1^{\text{\'et}}(X))$ is not itself an isomorphism, it is injective with dense image, if we give these two sets appropriate topologies. Take the point of view that the difference between $\operatorname{Mor}(Y,X) \rightarrow \operatorname{Hom}(\pi_1(Y),\pi_1(X))$ being an isomorphism and being injective with dense image is ``not very important,'' i.e. not suggestive of the presence of another sort of fundamental group. 

On the other hand, the local conditions inherent in the definition of the Selmer group are perhaps the result of some other sort of fundamental group. More explicitly, note that the rational points on an anabelian scheme are conjectured to be in bijection with $\Hom(\Gal(\kbar/k), \pi_1^{\text{\'et}})$ (Grothendieck's Section Conjecture).  However, for an elliptic curve, conditions must be imposed on an element of $\Hom(\Gal(\kbar/k), \pi_1^{\text{\'et}})$ for the element to come from a rational point. Explicitly, let $E$ be an elliptic curve over $k$, and let $S^n(E/k)$ and $\Sha(E/k)$ be the $n$-Selmer group and Shafarevich-Tate group of $E/k$ respectively. The exact sequence $$0 \rightarrow E(k)/n E(k) \rightarrow  S^n(E/k) \rightarrow \Sha(E/k)[n] \rightarrow 0$$ gives the exact sequence

$$
0 \rightarrow \varprojlim_n E(k)/n E(k) \rightarrow \varprojlim_n S^n(E/k) \rightarrow \varprojlim_n \Sha(E/k)[n] \rightarrow 0.
$$

Thus if $\Sha(E/K)$ has no non-zero divisible elements,  

\begin{equation}\label{completed_rational_pts_equals_completed_Selmer_group}
\varprojlim_n E(k)/n E(k) \cong \varprojlim_n S^n(E/k). 
\end{equation}

It is not hard to see that $\Hom(\Gal(\kbar/k), \pi_1^{\text{\'et}}) \cong H^1(\Gal(\kbar/k), \varprojlim_n E[n])$ and that  $\varprojlim_n S^n(E/k)$ is naturally a subset of $H^1(\Gal(\kbar/k), \varprojlim_n E[n])$. Think of $\varprojlim_n S^n(E/k)$ as a subset of $\Hom(\Gal(\kbar/k), \pi_1^{\text{\'et}})$ cut out by local conditions, as in the definition of the Selmer group. Any rational point of $E$ must be in this subset.

 Furthermore, if $\Sha(E/K)$ has no non-zero divisible elements, equation (\ref{completed_rational_pts_equals_completed_Selmer_group}) can be interpreted as saying that the (profinite completion of the) rational points of $E$ are the elements of $\Hom(\Gal(\kbar/k), \pi_1^{\text{\'et}})$ satisfying these local conditions.
 
 We ask if it is really necessary to exclude elliptic curves from the algebraic $K(\pi,1)$'s,  or if there is another sort of covering space, another sort of ``loop up to homotopy,'' producing a fundamental group which does characterize elliptic curves. For instance, if $\Sha(E/K)$ has no non-zero divisible elements, this example suggests that this new sort of fundamental group only needs to produce local conditions, perhaps by considering some sort of localization of the elliptic curve.


\begin{thebibliography}{[EGA III1]}


\bibitem[Ar]{artin} M. Artin, {\em Grothendieck Topologies}. 
Lecture Notes, Harvard University Math. Dept., Cambridge, Mass 1962.

\bibitem[AM]{am} M. Artin and B. Mazur, {\em Etale Homotopy},
Lecture Notes in Math.\ {\bf 100}, Springer-Verlag, Berlin-New York 1969.

\bibitem[AT]{artintate} E. Artin and J. Tate, {\em Class Field Theory},
W. A. Benjamin, Inc., New York-Amsterdam, 1968.

\bibitem[Bi1]{biss1} D. K. Biss, {\em A generalized approach to the fundamental
group}, Amer.\ Math.\ Monthly {\bf 107} (2000), no.\ 8, 711--720.

\bibitem[Bi2]{biss2} D. K. Biss, {\em The topological fundamental
group and generalized covering spaces}, Topology Appl.\ {\bf 124} (2002), no.\ 3, 355--371.

\bibitem[Bo]{boggi} M. Boggi, {\em Profinite Teichm\"uller theory},
Math.\ Nachr.\ {\bf 279} (2006), no.\ 9--10, 953--987.

\bibitem[BLR]{neronmodels}  S. Bosch, W. Lutkebohmert,  M. Raynaud, {\em N\'eron Models}, Springer, 1990.

\bibitem[BT]{bt} F. Bogomolov and Y. Tschinkel, {\em Unramified correspondences},
 in {\em Algebraic number theory and algebraic geometry}, 17--25, 
Contemp.\ Math., 300, Amer.\ Math.\ Soc., Providence, RI, 2002.

\bibitem[D]{deligne} P. Deligne, {\em Le groupe fondamental de la droite projective moins trois points}, in {\em Galois groups over $\Q$ (Berkeley, CA, 1987)}, 79--297, 
Math. Sci. Res. Inst. Publ., 16, Springer, New York, 1989. 

\bibitem[SGA 3]{sga3} M. Demazure and A. Grothendieck, {\em S\'eminaire de G\'eom\'etrie Alg\'ebrique du Bois Marie 1962/64, Sch\'emas en Groupes II,} Lecture Notes in Mathematics 152, Springer-Verlag, 1970. 

\bibitem[E]{ehrenpreis} L. Ehrenpreis, {\em Cohomology with bounds},
in {\em Symposia Mathematica, Vol.\ IV (INDAM, Rome, 1968/69)} pp.\ 389--395, Academic Press, London.

\bibitem[EH]{eh} H. Esnault and P. H. Hai, {\em The fundamental groupoid scheme and applications}, preprint 2006, arXiv:math/0611115v2.

\bibitem[Ge]{Geyer} W.-D. Geyer,  {\em Unendliche algebraische Zahlk\"orper, \"uber denen jede Gleichung aufl\"osbar von beschr\"ankter Stufe ist}, J. Number Theory {\bf 1} (1969), 346--374.

\bibitem[Gi]{Gille} P. Gille, {\em Le groupe fondamental sauvage d'une courbe affine en caract\'etristique $p>0$} in {\em Courbes semi-stables et groupe fondamental en g\'eom\'etrie alg\'ebrique (Luminy, 1998)}, 217--231, Progr.\ Math.\ {\bf 187}, Birkh\"auser, Basel, 2000.

\bibitem[EGA II]{ega2} A. Grothendieck, 
{\em \'El\'ements de g\'eom\'etrie
alg\'ebrique II. \'Etude globale \'el\'ementaire de quelques classes de morphismes}, IHES Publ.\ Math.\ No.\ 11, 1961.

\bibitem[EGA III${}_1$]{ega3} A. Grothendieck, 
{\em \'El\'ements de g\'eom\'etrie
alg\'ebrique III,  Premi\`ere Partie}, IHES Publ.\ Math.\ No.\ 11, 1961.

\bibitem[EGA IV${}_2$]{ega4ii} A. Grothendieck, {\em \'El\'ements de g\'eom\'etrie
alg\'ebrique IV:  \'Etude locale des sch\'emas et des morphismes
de sch\'emas, Seconde partie}, IHES Publ.\ Math.\ No.\ 24, 1965.

\bibitem[EGA IV${}_3$]{ega4iii} A. Grothendieck, {\em \'El\'ements de g\'eom\'etrie
alg\'ebrique IV:  \'Etude locale des sch\'emas et des morphismes
de sch\'emas, Troisi\`eme partie}, IHES Publ.\ Math.\ No.\ 28, 1967.

\bibitem[EGA IV${}_4$]{ega4iv} A. Grothendieck, {\em \'El\'ements de g\'eom\'etrie
alg\'ebrique IV:  \'Etude locale des sch\'emas et des morphismes
de sch\'emas, Quatri\`eme partie}, IHES Publ.\ Math.\ No.\ 32, 1967.

\bibitem[Gr1]{Groth_Brief_Faltings} A. Grothendieck, {\em Brief an G. Faltings.} (German) [Letter to G. Faltings] 
With an English translation on pp. 285--293. London Math. Soc. Lecture Note Ser., 242, Geometric Galois actions, 1, 49--58, Cambridge Univ. Press, Cambridge, 1997.

\bibitem[Gr2]{Esquisse} A. Grothendieck, {\em Esquisse d'un programme.} (French. French summary) [Sketch of a program] 
With an English translation on pp. 243--283. London Math. Soc. Lecture Note Ser., 242, Geometric Galois actions, 1, 5--48, Cambridge Univ. Press, Cambridge, 1997.

\bibitem[SGA1]{sga1} A. Grothendieck (dir.), {\em Rev\^etements \'etales et groupe fondamental (SGA1)}, Documents Math\'ematiques 3, Soc.\ Math.\ Fr., Paris, 2003.

\bibitem[KM]{km} J. Kahn and V. Markovic, {\em Random ideal triangulations
and the Weil-Petersson distance between finite degree covers of punctured
Riemann surfaces}, preprint 2008, arXiv:0806.2304v1.

\bibitem[Le]{lenstra} H. W. Lenstra, {\em Galois Theory for Schemes},
course notes available from the server of the Universiteit Leiden Mathematics Department, http://websites.math.leidenuniv.nl/algebra/GSchemes.pdf, Electronic third edition: 2008.

\bibitem[MS]{ms} V. Markovic and D. \v Sari\'c,
{\em Teichm\"uller mapping class group of the universal
hyperbolic solenoid}, Trans.\ Amer.\ Math.\ Soc.\ {\bf 358} (2006), no.\ 6,
2637--2650.

\bibitem[Mc]{mcmullen} C. McMullen, {\em Thermodynamics, dimension,
    and the Weil-Petersson metric}, Invent.\ Math.\ {\bf 173} (2008), no.\ 2, 365--425.



\bibitem[Mi]{milnebook} J. S. Milne, {\em \'Etale cohomology},
  Princeton Math.\ Ser., 33, Princeton U. P., Princeton, N.J., 1980.

\bibitem[Moc]{mochizuki} S. Mochizuki, {\em Correspondences on hyperbolic curves},
preprint, available at \verb+http://www.kurims.kyoto-u.ac.jp/~motizuki/papers-english.html+.

\bibitem[Mor]{Morel} F. Morel, {\em $\mathbb{A}^1$-algebraic topology over a field}, preprint, available at \verb+http://www.mathematik.uni-muenchen.de/~morel/preprint.html+. 

\bibitem[MM]{MM} 
J. W. Morgan and I. Morrison, {\em A van Kampen theorem for weak joins},  Proc.\ London Math.\ Soc.\ (3)  {\bf 53}  (1986),  no.\ 3, 562--576.


\bibitem[Mur]{murre}
J. P. Murre, {\em Lectures on an Introduction to Grothendieck's theory of the 
Fundamental Group}, notes by S. Anantharaman, TIFR Lect.\ on Math., no.\ 40,
TIFR, Bombay, 1967.

\bibitem[NSW]{coh_num_fields} J. Neukirch, A. Schmidt and K. Wingberg, {\em Cohomology of number fields}, Second edition. Grundlehren der Mathematischen Wissenschaften, 323. Springer-Verlag, Berlin, 2008. 

\bibitem[N1]{nori_compositio} {M. Nori}, {\em On the
     representations of the fundamental group}, {Compositio Math.},
   {\bf 33}, (1976), no.\ {1}, {29--41}, {\MR{0417179 (54 \#5237)}}.


\bibitem[N2]{nori} M. Nori, {\em The fundamental group-scheme}, 
Proc.\ Indian Acad.\ Sci.\ (Math.\ Sci.) {\bf 91}, no.\ 2, 1982, 73--122.

\bibitem[Oo]{oort} F. Oort,  {\em The algebraic fundamental group}, in {\em Geometric 
Galois Actions}, 1, 67--83, London Math.\ Soc.\ Lecture Note Ser., 242, Cambridge UP, Cambridge, 1997.

\bibitem[Pa]{pardini} R. Pardini, {\em Abelian covers of algebraic varieties},
J.\ Reine Angew.\ Math.\ {\bf 417} (1991), 191--213.

\bibitem[Po]{pop} F. Pop, {\em Anabelian Phenomena in Geometry and Arithmetic}, 2005 notes, available at 
\verb+http://modular.math.washington.edu/swc/aws/notes/files/05PopNotes.pdf+.

\bibitem[PSS]{Shaf} I. I. Piatetski-Shapiro and I. R. Shafarevich,
{\em Galois theory of transcendental extensions and uniformization}, in {\em Igor R. Shafarevich:  Collected Mathematical Papers}, Springer-Verlag, New York, 1980, 387-421.


\bibitem[RZ]{RZ} L. Ribes and P. Zalesskii, {\em Profinite Groups}, Ergebnisse der Mathematik und ihrer Grenzgebiete, 40, Springer, Berlin, 2000.

\bibitem[Se1]{serreRE} J.-P. Serre, {\em Construction de rev\^etements \'etales de la droite affine en caract\'eristique $p$.}
C. R. Acad. Sci. Paris S\'er. I Math. 311 (1990), no. 6, 341--346.

\bibitem[Se2]{serreGP} J.-P. Serre, {\em Groupes proalg\'ebriques},
IHES Publ.\ Math.\ No.\ 7, {\bf 1960}.

\bibitem[Se3]{Serre} J.-P. Serre, {\em Galois cohomology}, P. Ion trans.,
Springer-Verlag, Berlin, 2002.

\bibitem[Sp]{spanier} E. Spanier, {\em Algebraic Topology}, McGraw-Hill Book Co., New York, 1966.

\bibitem[Su]{sullivan} D. Sullivan, {\em Linking the universalities of
Milnor-Thurston, Feigenbaum, and Ahlfors-Bers}, in {\em Topological
Methods in Modern Mathematics (Stony Brook, NY, 1991)}, 543--564, Publish
or Perish, Houston TX, 1993.

\bibitem[Sz1]{SzLuminy} T. Szamuely {\em Le th\'eor\`eme de Tamagawa I}, in {\em Courbes semi-stables et groupe fondamental en g\'eom\'etrie alg\'ebrique (Luminy, 1998)}, 185--201, Progr.\ Math.\ {\bf 187}, Birkh\"auser, Basel, 2000.

\bibitem[Sz2]{sz} T. Szamuely, {\em Galois Groups and Fundamental Groups}, Cambridge Studies in Advanced Mathematics, vol. 117, to be published by Cambridge University Press in 2009.

\end{thebibliography}
\end{document}